\documentclass{IEEEtran4PSCC}
\usepackage{scrextend}
\usepackage[utf8]{inputenc}
\usepackage{bbm}
\usepackage{bm}
\usepackage{xcolor}
\usepackage{multirow}
\usepackage{bbold}
\usepackage{siunitx}
\usepackage{mathtools}
\usepackage{enumitem}
\usepackage{algorithm}
\usepackage{algpseudocode}

\newcommand{\pap}[1]{{\leavevmode\color{black}#1}}

\usepackage{cite}
\renewcommand{\baselinestretch}{0.9731}

\let\Right\right
\let\Left\left
\makeatletter
\def\right#1{\Right#1\@ifnextchar){\!\right}{}}
\def\left#1{\Left#1\@ifnextchar({\!\left}{}}
\makeatother

\ifCLASSOPTIONcompsoc
    \usepackage[caption=false, font=normalsize, labelfont=sf, textfont=sf]{subfig}
\else
\usepackage[caption=false, font=footnotesize]{subfig}
\fi

\ifCLASSINFOpdf
   \usepackage[pdftex]{graphicx}
\else
   \usepackage[dvips]{graphicx}
\fi
\usepackage{multirow}
\makeatletter
\let\old@ps@headings\ps@headings
\let\old@ps@IEEEtitlepagestyle\ps@IEEEtitlepagestyle
\def\psccfooter#1{%
    \def\ps@headings{%
        \old@ps@headings%
        \def\@oddfoot{\strut\hfill#1\hfill\strut}%
        \def\@evenfoot{\strut\hfill#1\hfill\strut}%
    }%
    \def\ps@IEEEtitlepagestyle{%
        \old@ps@IEEEtitlepagestyle%
        \def\@oddfoot{\strut\hfill#1\hfill\strut}%
        \def\@evenfoot{\strut\hfill#1\hfill\strut}%
    }%
    \ps@headings%
}
\makeatother

\begin{document}

\title{Adaptive Power Flow Approximations with Second-Order Sensitivity Insights}

\author{%
  \IEEEauthorblockN{%
    Paprapee Buason$^{\dagger}$\textsuperscript{\textsection},
    Sidhant Misra$^{\dagger}$,
    Jean-Paul Watson$^{\ast}$, and
    Daniel K. Molzahn\textsuperscript{\textsection}%
  }%
}

\maketitle

\thanksto{\noindent
$^{\dagger}$: Los Alamos National Laboratory. \{buason, sidhant\}@lanl.gov. They are supported by the Advanced Grid Modeling Program of the Office of Electricity of the U.S. Department of Energy.\\
\textsection: School of Electrical and Computer Engineering, Georgia Institute of Technology. molzahn@gatech.edu. They were supported by the National Science Foundation under grant number 2023140.\\
$^{\ast}$: Lawrence Livermore National Laboratory. watson61@llnl.gov. This work was performed under the auspices of the U.S. Department of Energy by Lawrence Livermore National Laboratory under Contract DE-AC52-07NA27344 and was supported in by the Advanced Grid Modeling Program of the Office of Electricity of the U.S. Department of Energy.}

\begin{abstract}

\textcolor{black}{The power flow equations are fundamental to power system planning, analysis, and control. However, the inherent non-linearity and non-convexity of these equations present formidable obstacles in problem-solving processes. To mitigate these challenges, recent research has proposed adaptive power flow linearizations that aim to achieve accuracy over wide operating ranges. The accuracy of these approximations inherently depends on the curvature of the power flow equations within these ranges, which necessitates considering second-order sensitivities. In this paper, we leverage second-order sensitivities to both analyze and improve power flow approximations. We evaluate the curvature across broad operational ranges and subsequently utilize this information to inform the computation of various sample-based power flow approximation techniques. Additionally, we leverage second-order sensitivities to guide the development of rational approximations that yield linear constraints in optimization problems. This approach is extended to enhance accuracy beyond the limitations of linear functions across varied operational scenarios.}

\end{abstract}

\section{Introduction} \label{sec:introduction}

\textcolor{black}{Accurately modeling the steady-state relationships among the bus voltages, line power flows, and active and reactive power injections, the AC power flow (AC-PF) equations are central to power system analyses. As such, these equations are essential for network planning and decision-making processes. However, the inherent non-linearity of the AC-PF equations presents formidable challenges when embedded in various optimization problems. Addressing complex issues such as uncertainty quantification for managing stochastic loads or tackling mixed-integer and bilevel problems for unit commitment, expansion planning, sensor placement, etc. is greatly impeded by the non-convex nature of the AC-PF equations.}

To obtain tractable formulations, a variety of power flow approximations are often used in place of the AC-PF equations within optimization problems. The DC approximation for high-voltage transmission systems \cite{stott2009}, LinDistFlow for distribution systems \cite{baran1989}, and the first-order Taylor expansion are among the most commonly used approximations. \textcolor{black}{Beyond these, there is a rich literature on alternative approximation approaches and methods to enhance existing approximations. Recent work includes, for instance, methods for improving DC power flow approximations by optimizing bias and coefficient parameters~\cite{taheri_molzahn-dcparam} and incorporating voltage magnitude and reactive power considerations~\cite{yang2018}. Other recent work uses an optimization strategy to select the most suitable choice of voltage representation~\cite{fan2021, BUASON2022}. Reference~\cite{fnt} surveys recent developments in power flow approximations.}

Unsurprisingly, many power flow approximations are linear in order to obtain linear constraints that significantly improve tractability of optimization formulations relative to the nonlinear AC-PF equations. For complex and large-scale mixed-integer and multi-level problems, preserving linearity may be the only viable strategy to obtain tractable solution methods. 

Many well-known power flow approximations rely on general assumptions about operating ranges and system characteristics. For example, the DC approximation relies on negligible voltage variations and losses. The generality of these approximations, however, can come with the tradeoff of large errors that lead to constraint violations and sub-optimal cost. To address this issue, several recent works have considered \emph{adaptive} power flow approximations that are tailored to be accurate for a particular system over a specified operating range. These include optimization-based approaches~\cite{misra2018optimal,muhlpfordt2019optimal,taheri_molzahn-dcparam,taheri_gupta_molzahn-opt_lindistflow} and data-driven approaches~\cite{BUASON2022,buason2023datadriven, buason_misra_molzahn-cbla, CHEN2022108573, fan2021, liu2018}; see~\cite{10202779, jia2023tutorial1, jia2023tutorial2} for recent~surveys. 

\textcolor{black}{These adaptive power flow approximations invest computing time to calculate linearization coefficients that are then used to better approximate nonlinear optimization problems by linear problems. In contrast, traditional power flow formulations like DC power flow, LinDistFlow, and first-order Taylor expansion do not require up-front time for computing the linearization coefficients, but may result in lower accuracy when used to linearly approximate non-linearities in optimization problems~\cite{baker2021, dvijotham_molzahn-cdc2016}. In other words, there is an inherent computational tradeoff for adaptive power flow approximations between up-front computing time and accuracy of a linearized optimization problem. This tradeoff is most favorable in settings characterized by an offline/online split, where ample time is available offline for calculating the coefficients in order to perform fast online computations (e.g., employing a day-ahead forecast of load demands in offline computations of linearization coefficients that are then used online in real-time dispatch computations). Moreover, these methods can address optimization problems that would otherwise be intractable or exceedingly challenging with nonlinear AC power flow models, such as AC Unit Commitment~\cite{castillo2016}, resilient infrastructure planning with AC power flow models~\cite{bhusal2020, austgen2023comparisons, haag2024, owen_aquino_talkington_molzahn-EVacuation_feeder}, and optimization of sensor placement~\cite{buason2023datadriven}. In these settings, the up-front investment of computational effort to calculate the linearization coefficients can be worthwhile since incorporating the nonlinear AC power flow equations poses severe challenges for these problems and typical linearizations can result in significant inaccuracies.}

To further improve the utility of adaptive linear approximations for optimization problems, \cite{BUASON2022} introduced the concept of \emph{conservative} linear approximations (CLAs) that are designed to either under- or over-estimate a quantity of interest such as bus voltage magnitude. Since the goal of these approximations is to formulate inequality constraints in an optimization problem (e.g., satisfying limits on voltage magnitudes), introducing conservativeness during their construction yields improved performance in terms of constraint enforcement.

Data-driven formulations, which rely on randomly sampled points from a given operating range, are the most widely used methods for constructing adaptive power flow approximations. Tractability, scalability, and suitability for parallel implementation are the primary reasons for their popularity. Naturally, the accuracy of the approximations heavily depends on the quantity and quality of samples used in their construction.

This paper develops an \emph{importance sampling} method to construct linear and conservative linear approximations. The method aims to include the most informative samples to quickly improve approximation quality while keeping the training process efficient. Our approach is based on second-order sensitivity information. By preferentially including more samples from a relatively low-dimensional subspace with high curvature, we obtain highly accurate linear approximations using far fewer samples compared to random sampling. 

\textcolor{black}{Second-order derivatives are key to informing the importance sampling method. This paper is, to the best of our knowledge, the first to present the second-order derivatives of the voltage magnitudes with respect to active and reactive power injections. Other literature considers the second-order derivatives of the power injections with respect to voltage magnitudes and angles~\cite{matpowerTechNote2}, the second-order derivatives of network losses and incremental costs with respect to active and reactive power injections~\cite{lee2018}, the second-order derivatives of network losses with respect to voltage magnitudes and angles~\cite{martin2016}, and the second-order derivatives of the voltage phasors with respect to the complex power injections and their conjugates at the no-load operating point~\cite{jabr2019}. These all differ from the sensitivities we derive and are unsuitable for the importance sampling method developed in this paper.}

We also introduce a class of approximations based on rational functions with linear numerator and denominator. The main motivation behind using these rational functions is that they result in linear inequality constraints when used in an optimization formulation. At the same time, expanding the approximating class of functions from linear to rational provides better approximation error. We first introduce a multivariate generalization of the Pad\'e approximant which can be constructed at a given operating point using the second-order sensitivities and serves as a rational generalization of the second-order Taylor expansion~\cite{Cuyt1984}. \textcolor{black}{As far as we are aware, this is the first application of Pad\'e approximants to power flow approximations.}\footnote{\textcolor{black}{We note that Pad\'e approximants are used to improve the convergence of holomorphic embedding methods for solving power flow problems~\cite{trias2018helm}, but this differs from both the motivation and mathematical context of this paper.}} We then provide a data-driven training method similar to \cite{BUASON2022} to construct rational and conservative rational approximations over an operating range.

\textcolor{black}{In summary, the main contributions of this paper are:
\begin{enumerate}[wide, labelwidth = 0pt, labelindent = 0pt, label=(\textit{\roman*})]
    \item Derivation of the second-order sensitivity matrix of voltage magnitudes with respect to complex power injections.
    \item Analysis of second-order sensitivities to draw informative samples for adaptive power flow approximations.
    \item Introduction of a multivariate generalization of the Pad\'e approximant that utilizes second-order sensitivities to enhance power flow approximations formulated as rational functions.
    \item Development of an adaptive approach to constructing conservative rational approximations using a constrained regression formulation.
    \item Numerical analysis of the multivariate  generalization of the Pad\'e approximant and the rational approximations across various test cases, demonstrating their application to the optimal power flow problem.
\end{enumerate}
}

The paper is organized as follows: Section~\ref{sec:background} covers background material on the power flow equations as well as adaptive linear and conservative linear approximations. Section~\ref{sec:adaptive sampling} describes the importance sampling approach for data-efficient construction of linear approximations. Section~\ref{sec:linear approx} introduces rational approximations, including the multivariate generalization of the Pad\'e approximant for the power flow equations. Section~\ref{sec:simulation} provides numerical results demonstrating the improvements obtained from our approach. Section~\ref{sec:future work} concludes the paper along with directions for our future work.

\section{Background} \label{sec:background}

In this section, we first present the power flow equations and then describe a so-called ``conservative linear approximation'' (CLA) approach for both voltage magnitudes, as we previously proposed, and a new extension to current flows.

\subsection{The power flow equations}\label{sub:pf}

Consider an $N$-bus power system. Let the length-$N$ vectors $P$, $Q$, and $V \angle \theta$ denote the active and reactive power injections and voltage phasors, respectively, at each bus. Designate a reference bus where the voltage phasor is set to $1 \angle 0$\textdegree~per unit. 
We use the subscript $(\cdot)_{i}$ to represent a quantity at bus $i$ and the subscript $(\cdot)_{ik}$ to represent a quantity from or connecting bus $i$ to $k$. The AC power flow equations for bus~$i$ are:
\begin{subequations}
\label{eq:power_flow}
\vspace{-0.25em}
\begin{align}
	P_i &= V_i^2 G_{ii} + \sum_{k \in \mathcal B_i} V_i V_k(G_{ik}\cos \theta_{ik} + B_{ik}\sin\theta_{ik}), \label{eq:dP} \\
	Q_i &= -V_i^2 B_{ii} + \sum_{k \in \mathcal B_i} V_i V_k(G_{ik}\sin\theta_{ik} - B_{ik}\cos\theta_{ik}), \label{eq:dQ}
\end{align}
\end{subequations}
where $G +jB$ (with $j=\sqrt{-1}$) is the admittance matrix.


\subsection{Conservative linear approximations} \label{sub:cla}
\begin{figure}[t]
	\centering 
	\includegraphics[trim=0.5cm 0.5cm 0.2cm 0.8cm, clip, width=0.8\linewidth]{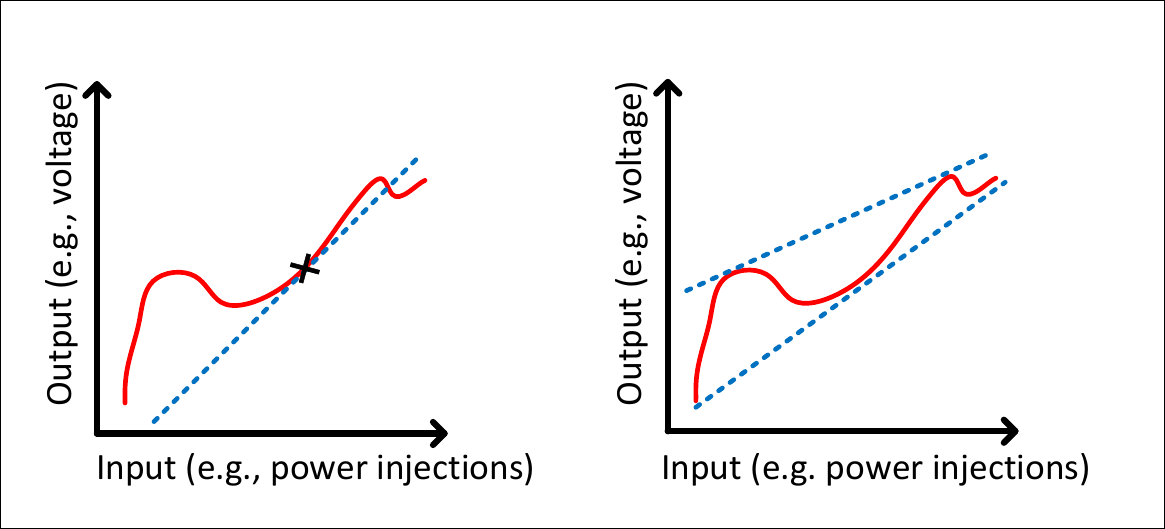} 
	\caption{A conceptual example of a traditional linear approximation (left) and a conservative linear approximation (right). The solid line represents the nonlinear function of interest. The dotted line in the left figure is a traditional first-order Taylor approximation around point $\times$ while the dotted top (bottom) line in the right figure is an over- (under-) estimating approximation.}
	\label{fig:cla_fig}
 \vspace{-1.5em}
\end{figure}

To model voltage and current limits, many optimization problems include the power flow equations~\eqref{eq:power_flow} as constraints. The non-linearity of these equations contributes to the complexity of the problems, often making them difficult to solve. To address this challenge, we previously introduced a sample-based CLA approach in~\cite{BUASON2022}. This linear approximation seeks to over- or under-estimate a specified quantity of interest (see Fig.~\ref{fig:cla_fig}). CLAs approximate the nonlinear power flow equations using a sample-based approach that enables parallel computation, i.e., the CLA of each quantity of interest can be computed concurrently. Constructing a CLA involves sampling power injections over an operating range of interest, computing the power flow equations for each sample, and then solving a constrained-regression problem. This approach tailors the approximation to a particular operating range and system of interest. Moreover, considering more complex components, such as tap-changing transformers and smart inverters, is straightforward since their behavior can be incorporated into the sampled power flow solutions~\cite{buason2023datadriven}. CLAs thus yield simplified optimization problems that are suitable for commercial optimization solvers. Finally, when applied in an optimization context, CLAs have a key advantage: one may ensure satisfaction of nonlinear constraints \emph{while only enforcing linear inequalities} (assuming that the CLAs are indeed conservative). In this section, we next revisit the CLA approach.

We denote vectors and matrices in bold. Consider some quantity of interest (e.g., the voltage magnitude at a particular bus or current flow on a certain line) that we generically denote as $\beta$. An \emph{overestimating} CLA is given by the linear expression:
\begin{equation}\label{eq:cla_template}
a_{0} + \bm{a}_{1}^T\begin{bmatrix}
\bm{P} \\
\bm{Q}
\end{bmatrix},
\end{equation}
constructed such that the following relationship is satisfied for power injections $\bm{P}$ and $\bm{Q}$ within a specified range:
\begin{equation}
\beta \leq a_{0} + \bm{a}_{1}^T\begin{bmatrix}
\bm{P} \\ \bm{Q}
\end{bmatrix},
\label{eq:cla_setup}
\end{equation}
where the superscript $T$ denotes the transpose. Assuming that~\eqref{eq:cla_setup} is indeed satisfied, one may ensure satisfying the constraint $\beta \leq \beta^{\text{max}}$, where $\beta$ is the output of a nonlinear function such as the implicit system of nonlinear AC power flow equations~\eqref{eq:power_flow} and $\beta^{\text{max}}$ is a specified upper bound, by instead enforcing the linear constraint $a_{0} + \bm{a}_{1}^T\begin{bmatrix}
\bm{P} \\
\bm{Q}
\end{bmatrix} \leq \beta^{\text{max}}$. 

To determine the coefficients of the affine function of power injections in~\eqref{eq:cla_template}, we solve the following regression problem:
\begin{subequations} \label{eq:regression}
\vspace{-0.35em}
\begin{align}
\label{eq:regression_loss}
&\min_{a_{0}, \ \bm{a}_{1}}  \quad \frac{1}{M}\sum_{m = 1}^M \mathcal{L}\left(\beta_{m} - \left(a_{0} + \bm{a}_{1}^T\begin{bmatrix}
\bm{P}_m \\ \bm{Q}_m
\end{bmatrix}\right)\right) \\
&\text{s.t.} \quad \beta_{m} - \left(a_{0} + \bm{a}_{1}^T\begin{bmatrix}
\bm{P}_m \\ \bm{Q}_m
\end{bmatrix}\right) \leq 0, \quad m = 1,\ldots, M,\label{eq:regression_over}
\end{align}%
\label{eq:cla}%
\end{subequations}
where the subscript $m$ denotes the $m^{\text{th}}$ sample, $M$ is the number of samples, and $\mathcal{L}(\,\cdot\,)$ represents a loss function, e.g., absolute value for the $\ell_1$ loss and the square for the squared-$\ell_2$ loss. In this paper, our quantities of interest ($\beta$) are the magnitudes of voltages ($V$) and current flows ($I$).  

Underestimating CLAs are constructed in the same manner as~\eqref{eq:regression} except that the direction of the inequality in~\eqref{eq:regression_over} is flipped. One could also compute a linear approximation (LA) that is not conservative (minimizes approximation errors without consistently under- or over-estimating) by solving~\eqref{eq:regression} without enforcing~\eqref{eq:regression_over}.

\section{Importance sampling method\\ for power flow approximations} \label{sec:adaptive sampling}

The accuracy and conservativeness of adaptive power flow approximations like CLA depend on the set of samples, with more samples leading to higher accuracy and increased confidence in conservativeness. However, using many samples can be computationally challenging, which is problematic when power flow approximations need to be recomputed frequently during rapidly changing operating conditions.

Thus, approximation methods such as CLA can significantly benefit from the concept of importance sampling.\footnote{\textcolor{black}{This section describes an importance sampling method in terms of the conservative linear approximations introduced in the prior section. This method also applicable to the conservative rational approximations presented later in Section~\ref{sec:linear approx}.}} The samples that provide the most valuable information for improving the quality of a linearization are often associated with regions of the power flow manifold that exhibit the most curvature. We propose an importance sampling method based on the second-order sensitivities of the power flow equations. This sampling method can achieve high accuracy with fewer samples compared to naive random sampling approaches. 

For brevity, we rewrite the power flow equations in~\eqref{eq:power_flow} as:
\begin{align}
    \bm{x} = g(\bm{y}) \label{eq:pf_simplified},    
\end{align}
where $g$ is a vector-valued function and
\begin{equation}\label{eq:x_y_def}
    \bm{x} = \begin{bmatrix}\bm{P} \\ \bm{Q} \end{bmatrix}, \quad     \bm{y} = \begin{bmatrix} \bm{\theta} \\ \bm{V}\end{bmatrix}.
\end{equation}

The second-order sensitivity matrix for a specific quantity of interest $y_k$, denoted as $\bm{\Lambda}_{y_k}$, takes the following form:
\begin{equation}
\label{eq:Lambda}
\bm{\Lambda}_{y_k} = 
    \begin{bmatrix}
        \cfrac{\partial^2 y_k}{\partial x_1 \partial x_1} &\cfrac{\partial^2 y_k}{\partial x_1 \partial x_2} &\cdots &\cfrac{\partial^2 y_k}{\partial x_1 \partial x_{2N}} \\
        \vdots &\vdots &\ddots &\vdots \\
        \cfrac{\partial^2 y_k}{\partial x_{2N} \partial x_1} &\cfrac{\partial^2 y_k}{\partial x_{2N} \partial x_2} &\cdots &\cfrac{\partial^2 y_k}{\partial x_{2N} \partial x_{2N}} 
    \end{bmatrix}.
\end{equation}
The appendix provides an explicit form for the second-order sensitivity matrix $\bm{\Lambda}_{y_k}$ in~\eqref{eq:Lambda}.

After obtaining the second-order sensitivity matrix, we analyze the power flow manifold's curvature by computing the singular value decomposition (SVD) of $\bm{\Lambda}_{y_k}$. A singular vector associated with a largest singular value indicates the direction of highest curvature. As we will demonstrate empirically in Section~\ref{sec:simulation}, the sensitivity matrices $\bm{\Lambda}_{y_k}$ typically have only a few significant singular values (i.e., these matrices are approximately low rank). Thus, the subspace spanned by only a few singular vectors gives a good characterization for the directions of highest curvature and hence suggests promising directions for prioritized sampling. By selecting a larger fraction of samples in the span of the dominant singular vectors, we seek to better inform the linear approximation to fit the most nonlinear components of the function. The second-order sensitivity analyses can also identify the convexity or concavity of the power flow equations at a specified operating point that can be used to further improve the sample selection.

\begin{figure}[t]
	\centering 
	\includegraphics[trim=0.5cm 0.2cm 0.4cm 0.3cm, clip, width=0.65\linewidth]{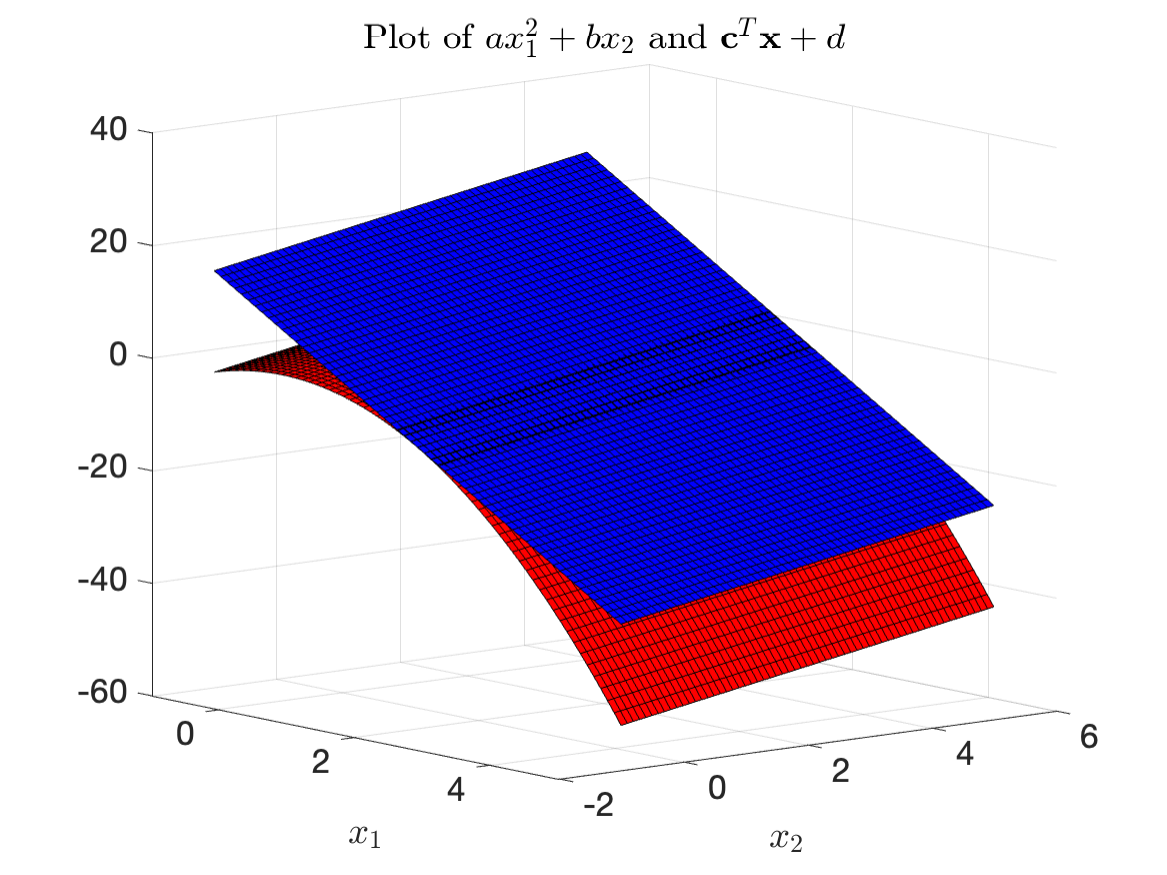} 
	\caption{An example of an overestimating linear function $\mathbf{c}^T\mathbf{x} + d$ (blue plane) of a quadratic function $ax_1^2 + bx_2$ (red manifold) where $\mathbf{x}^T = [x_1, \ x_2], a = -2, b = 2, \mathbf{c}^T = [-8, \ 2]$, and $d = 8$.}
	\label{fig:overestimate_plane_example}
	\vspace{-0.5em}
\end{figure}

The choice to draw samples along the span of singular vectors associated with the first few largest singular values of the second-order sensitivity matrix is conceptually illustrated via the simple three-dimensional plot in Fig.~\ref{fig:overestimate_plane_example}. This figure shows an example quadratic function that is nonlinear in the $x_1$ direction and linear in the $x_2$ direction along with an overestimating CLA. The CLA here minimizes the volume between the linear approximation and the nonlinear function while ensuring that the overestimating linear function consistently lies above the quadratic function. In a sample-based approach to computing the CLA, samples in the linear $x_2$ direction would provide little benefit while samples in the nonlinear $x_1$ direction would be informative. 

We present detailed results on power system test cases in Section~\ref{sec:simulation}, empirically showing that the second-order sensitivity matrices associated with voltage magnitudes are approximately low-rank and have all non-positive eigenvalues \textcolor{black}{indicating the concavity of voltage magnitudes as a function of power injections}. 
\textcolor{black}{This has significant implications for how to best select samples, especially for conservative power flow approximations. The conceptual example in Fig.~\ref{fig:overestimate_plane_example} helps illustrate this. The red surface has one direction with concave curvature ($x_1$) and one direction that is linear ($x_2$). To ensure conservativeness (i.e., the blue plane is always above the red surface), samples drawn near the middle of the range are most impactful; samples drawn near the corners of the operating range will be lower (and thus not affect the conservativeness requirement~\eqref{eq:regression_over}) due to the blue function's concavity. On the other hand, more extreme samples will have bigger contributions to the loss function in the regression problem~\eqref{eq:regression_loss}. Underestimating linear approximations have the converse behavior. Thus, our concavity results demonstrate that overestimating conservative linear approximations benefit most from samples near the middle of the operating range to ensure conservativeness with fewer extreme samples to maintain accuracy. 
Conversely, underestimating conservative linear approximations benefit from more extreme samples for conservativeness with fewer samples in the middle of the operating range for maintaining accuracy. In both cases, the low-rank nature provides a small number of directions where the power flow equations have the most substantial nonlinear behavior, thus providing the most impactful samples.}

\section{Rational Approximations of the\\ Power Flow Equations} \label{sec:linear approx}
The previously developed CLAs discussed in Section~\ref{sub:cla} are linear functions of the power injections. The primary goal of these linear functions is to maintain linearity in an optimization problem's constraints. However, there are more general classes of functions that are nonlinear but can be represented with linear constraints, specifically, rational functions with a linear numerator and a strictly positive linear denominator. With more degrees of freedom than linear functions, rational approximations can more accurately capture power flow nonlinearities. To exploit this more general class of functions, this section introduces the Pad\'e approximant and conservative rational approximations (CRAs) of the power flow equations. 

\subsection{Pad\'e approximant} \label{sub:Pade}
The Pad\'e approximant is a rational approximation that matches the first terms of a function's Taylor expansion~\cite{BAKER196121}. The Pad\'e approximant of a univariate function $f(x)$ is:
\begin{equation}
R_{x_0}(x) = \cfrac{a_0 + a_1 (x - x_0) + \cdots + a_m (x - x_0)^m}{1 + b_1 (x - x_0) + \cdots + b_n (x - x_0)^n}, \label{eq:Pade}
\end{equation}
where $m \geq 0$, $n \geq 1$, and both $m$ and $n$ are integers. Equation~$\eqref{eq:Pade}$ is called the $[m/n]$ Pad\'e approximant, which, \textcolor{black}{for univariate functions,} matches the $(m + n)^{\text{th}}$-degree Taylor series. To compute coefficients $a_0, \ldots, a_m$ and $b_1, \ldots, b_n$, we set $R_{x_0}(x)$ equal to the Taylor series and multiply both sides by the denominator (see~\cite{Cuyt1984}). By matching the coefficients of terms with degrees less than or equal to $m + n$, we can determine all the coefficients.

This paper focuses on the [1/1] Pad\'e approximant, which represents a ratio of linear functions. Despite being a nonlinear function, the Pad\'e approximant can be integrated into an optimization problem as a linear constraint. 
To be applicable to power systems settings, the input variable $\bm{x}$ must be a vector instead of a scalar. We therefore formulate a multivariate generalization of the Pad\'e approximant. Consider the constraint:
\begin{equation}
\label{eq:Pade_constraint1}
\cfrac{a_0 + \bm{a_1}^T \bm{x}}{1 + \bm{b_1} ^T \bm{x}} - U \leq 0,
\end{equation}
where $\bm{x}$ is a vector of decision variables and $U$ is an upper bound on the constrained quantity. In the multivariate generalization of the Pad\'e approximant shown in~\eqref{eq:Pade_constraint1}, $\bm{a_1}$ and $\bm{b_1}$ are vectors. For instance, we may seek to compute parameters $a_0, \bm{a_1}$, and $\bm{b_1}$ such that $(a_0 + \bm{a_1}^T \bm{x}) / (1 + \bm{b_1}^T\bm{x})$ approximates a voltage magnitude $V_i$ based on the power injections $\bm{x} = [\bm{P}^T\; \bm{Q}^T]^T$, with $U = V_i^\text{max}$ denoting the upper bound on the voltage at bus $i$. As long as $1 + \bm{b_1}^T \bm{x} > 0$, \eqref{eq:Pade_constraint1} can be reformulated as an equivalent linear constraint\footnote{\textcolor{black}{For the constraint~\eqref{eq:Pade_constraint1} around $\bm{x} = \bm{x_0}$, we can express an equivalent linear constraint as $(a_0 - U) + (\bm{a_1} - U \bm{b_1})^T \bm{x} \leq -(\bm{a_1} - U \bm{b_1})^T \bm{x_0}$.}}:
\begin{equation}
(a_0 - U) + (\bm{a_1} - U \bm{b_1})^T \bm{x} \leq 0 \label{eq:Pade_constraint2}.
\end{equation}

As we will show next, the additional degrees of freedom compared to a linear approximation enable the [1/1] multivariate Pad\'e approximant to better capture nonlinear behavior while still yielding a linear constraint.

The $[1/1]$ multivariate generalization of the Pad\'e approximant that we consider seeks to match the first- and second-order curvature of the second-order Taylor series approximation of the power flow equations.
Consider the following equation around $\bm{x} = \bm{x_0}$:
\begin{align}
\cfrac{a_0 + \bm{a_1}^T (\bm{x} - \bm{x_0})}{1 + \bm{b_1} ^T (\bm{x} - \bm{x_0})} &=  f(\bm{x_0}) + \nabla_{\bm{x}} f(\bm{x_0})^T (\bm{x} - \bm{x_0}) \notag \\ 
& \quad + \cfrac{1}{2}(\bm{x} - \bm{x_0})^T \bm{\Lambda}_{f}(\bm{x_0}) (\bm{x} - \bm{x_0}), \label{eq:Pade Taylor}
\end{align}
where $\nabla_{\bm{x}} f(\bm{x_0})$ is the gradient of $f$ at $\bm{x} = \bm{x_0}$ and $\bm{\Lambda}_{f}(\bm{x_0})$ is the second-order sensitivity matrix of $f$ at $\bm{x} = \bm{x_0}$. Multiplying by the denominator in~\eqref{eq:Pade Taylor} and comparing coefficients up to the quadratic term, we obtain:
\begin{subequations}
    \begin{align}
        a_0 &= f(\bm{x_0}), \\
        \bm{a_1} &= \nabla_{\bm{x}} f(\bm{x_0}) + f(\bm{x_0})\bm{b_1}, \\
        \bm{A} &= \bm{b_1} \nabla_{\bm{x}} f(\bm{x_0})^T + \cfrac{1}{2}\bm{\Lambda}_{f}(\bm{x_0}), \label{eq:second_order_Pade}
    \end{align}
\end{subequations}
for some skew-symmetric matrix $\bm{A}$.\footnote{If $\bm{A}$ is skew symmetric (i.e., $\bm{A} = -\bm{A}^T$), then $\bm{x}^T \bm{A} \bm{x}$ = 0 for any $\bm{x}$.} Since the matrix given by the outer product $\bm{b_1} \nabla_{\bm{x}} f(\bm{x_0})^T$ has rank 1, the matrix equation~\eqref{eq:second_order_Pade} is unsolvable, in general.\footnote{\textcolor{black}{In contrast to the univariate case where the [1/1] Pad\'e approximant matches the first- and second-order curvature of the second-order Taylor approximation, this unsolvability indicates that a [1/1] multivariate generalization of the Pad\'e approximant cannot replicate the first- and second-order curvature of the second-order Taylor expansion.}} Accordingly, there are many possible generalizations of univariate Pad\'e approximations to multivariate settings~\cite{cuyt1999}. 
We choose to approximately solve~\eqref{eq:second_order_Pade} in the Frobenius-norm sense. We first add the expression in~\eqref{eq:second_order_Pade} to  its transpose to obtain:
\begin{align}
   \bm{\underline{0}} &= \bm{b_1} \nabla_{\bm{x}} f(\bm{x_0})^T + \nabla_{\bm{x}} f(\bm{x_0}) \bm{b_1}^T  + \bm{\Lambda}_{f}(\bm{x_0}),  \label{eq:second_order_Pade_skew}
\end{align}
where $\bm{\underline{0}}$ is the all zeros matrix. The left-hand side equals $\bm{\underline{0}}$ since $\bm{A} + \bm{A}^T = \bm{0}$ for any skew-symmetric matrix $\bm{A}$. We then minimize the Frobenius norm (denoted as $\|\,\cdot\,\|_{\mathcal{F}}$) to select~$\bm{b_1}$:
\begin{equation}
    \bm{b_1} = \underset{\bm{b_1}}{\operatorname{arg\,min}} \;  \|\bm{b_1} \nabla_{\bm{x}} f(\bm{x_0})^T + \nabla_{\bm{x}} f(\bm{x_0}) \bm{b_1}^T  + \bm{\Lambda}_{f}(\bm{x_0})\|_{\mathcal{F}}. \label{eq:frobenius}
\end{equation}

We note that the second-order sensitivity $\bm{\Lambda}_{f}(\bm{x_0})$ discussed in Section~\ref{sec:adaptive sampling} is again needed here to compute the $[1/1]$ multivariate generalization of the Pad\'e approximant. By incorporating information from the second-order sensitivity, this approximant can better capture power flow nonlinearities while still yielding linear constraints in optimization problems. 

\subsection{Conservative rational approximations} \label{sub:cra}

The first-order Taylor series is the best-fitting linear function around a specified point. The CLAs in~\cite{BUASON2022} extend this concept to fit linear functions in an operating range of interest. Analogously, the [1/1] Pad\'e approximant is the best-fitting ratio of linear functions around a specified point. Motivated by the [1/1] multivariate generalization of the Pad\'e approximant, we next extend this concept to construct conservative rational approximations (CRAs) with linear numerator and denominator, called [1/1] CRAs, within an operating range of interest. 

A [1/1] CRA is defined as a rational function whose numerator and denominator are linear functions and, similar to CLAs, are conservative (under- or over-estimate the nonlinear power flow equations) in an operating range of interest. In an attempt to achieve this, we enforce conservativeness of the CRAs over a set of sampled power injections in this range. 

Analogous to~\eqref{eq:regression} for an overestimating CLA, the regression problem for computing an overestimating CRA is:
\begin{subequations} \label{eq:regression_socp}
\vspace{-0.35em}
\begin{align}
\min_{a_0,\bm{a_1},\bm{b_1}} & \quad  \frac{1}{M} \sum_{m = 1}^M \mathcal{L} \left(\beta_m - \cfrac{a_0 + \bm{a_1}^T \bm{x_m}}{1 + \bm{b_1}^T \bm{x_m}} \right) \label{eq:objective_socp}\\
\text{s.t.} \quad & a_0 + \bm{a_1}^T \bm{x_m} - \beta_m (1 + \bm{b_1}^T \bm{x_m}) \geq 0, \label{eq:conservative_socp}\\
\phantom{\text{s.t.}} \quad & 1 + \bm{b_1}^T \bm{x_m} \geq \epsilon, \label{eq:positive_denominator_socp}\quad m = 1,\ldots,M,
\end{align}
\end{subequations}
where $M$ again represents the number of samples.
The objective~\eqref{eq:objective_socp} minimizes the mismatch between the rational approximation and the value of the quantity of interest, again generically denoted as $\beta$, obtained from solving the power flow equations at each sampled set of power injections. 
%
Constraint~\eqref{eq:conservative_socp} ensures an overestimating property. Constraint~\eqref{eq:positive_denominator_socp}, where $\epsilon$ is a specified positive number, maintains the correct inequality sign in~\eqref{eq:conservative_socp}. The formulation for an underestimating CRA is given by~\eqref{eq:regression_socp} with a flipped inequality sign in~\eqref{eq:conservative_socp}. One could also compute a rational approximation (RA) that was not conservative by solving~\eqref{eq:regression_socp} without enforcing~\eqref{eq:conservative_socp}. Similar to the [1/1] multivariate generalization of the Pad\'e approximant, the rational approximations from solving~\eqref{eq:regression_socp} can better capture the curvature of the power flow manifold compared to linear approximations, here considering an operating range instead of a particular point.

Ideally, we would like to directly solve~\eqref{eq:regression_socp}. However, due to the objective \eqref{eq:objective_socp}, this is a nonlinear program that can be challenging to solve. Alternatively, we consider the following regression problem for computing an overestimating CRA with an $\ell_1$ loss function:
\begin{subequations} \label{eq:regression_linear}
\vspace{-0.35em}
\begin{align}
\min_{a_0,\bm{a_1},\bm{b_1}} \   & \frac{1}{M} \sum_{m = 1}^M w_m \left|\left(a_0 + \bm{a_1}^T \bm{x_m} - \beta_m (1 + \bm{b_1}^T \bm{x_m})\right)\right|\label{eq:regression_linear_objective}\\
\text{s.t.} \quad & a_0 + \bm{a_1}^T \bm{x_m} - \beta_m (1 + \bm{b_1}^T \bm{x_m}) \geq 0, \label{eq:conservative_linear}\\
& 1 + \bm{b_1}^T \bm{x_m} \geq \epsilon, \quad \label{eq:positive_denominator_linear} m = 1,\ldots,M,
\end{align}
\end{subequations}
%
%
where $w_m$ is a specified parameter vector that weights each term in the objective; \textcolor{black}{note that varying $w_m$ changes the solution to~\eqref{eq:regression_linear}}. If $w_m$ were equal to $(1 + \bm{b_1}^T \bm{x_m})^{-1}$, then~\eqref{eq:regression_linear} would be equivalent to~\eqref{eq:regression_socp} since multiplying by the weights $w_m$ would then compensate for multiplying by the denominators $1 + \bm{b_1}^T \bm{x_m}$ in~\eqref{eq:regression_linear_objective} relative to~\eqref{eq:objective_socp}. However, this choice of $w_m$ requires prior knowledge of $\bm{b_1}$, which is what we are trying to compute. 

We therefore employ an iterative algorithm that updates the weights $w_m$ at each iteration based on the values of $\bm{b_1}$ from the previous iteration. Let superscript $(\,\cdot\,)^k$ denote quantities at the $k^{\text{th}}$ iteration. This algorithm first initializes $\bm{w}$ as $\bm{w}^0$ (using, for instance, the $\bm{b_1}$ parameters computed from the [1/1] multivariate generalization of the Pad\'e approximant described in Section~\ref{sub:Pade}).  For each iteration~$k$, the algorithm then solves \eqref{eq:regression_linear} to compute $\bm{b_1}^k$ and updates $w_m^{k+1} = (1 + (\bm{b_1}^k)^T \bm{x_m})^{-1}$, repeating until $\|\bm{w}^k - \bm{w}^{k-1}\|_1$ reaches a specified tolerance. \textcolor{black}{In our numerical experience, this algorithm typically converges within five to ten iterations. Due to the scalability of linear programming solvers, this algorithm is much more tractable than directly solving~\eqref{eq:regression_socp}.}

\subsection{Computation processes, tractability, and parallelization} \label{sub:computation}

\begin{figure}[th!]
	\centering 
	\includegraphics[trim=0.5cm 0.9cm 0.5cm 0.5cm, clip, width=0.68\linewidth]{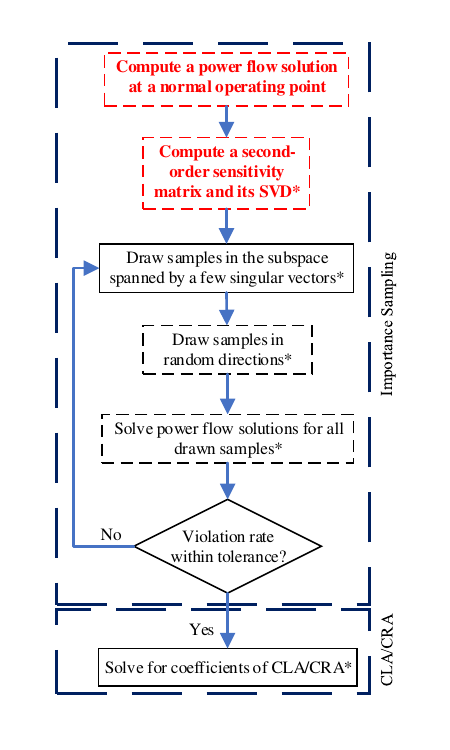} 
	\caption{\textcolor{black}{Flowchart depicting the computation processes for the importance sampling and CLA/CRA methods. Steps with $*$ are parallelizable. Red-dashed boxes highlight the computation of the second-order sensitivity matrix, while black-dashed boxes indicate sample drawing without importance sampling.}}
	\label{fig:algorithm}
 \vspace{-1.5em}
\end{figure}

\textcolor{black}{The computation processes for CRA include two key parts: selecting valuable samples using the importance sampling method and iteratively solving the regression problem in~\eqref{eq:regression_linear}. The detailed computational steps are illustrated in Fig.~\ref{fig:algorithm}.}

\pap{With the ability to parallelize all key computations, this paper's approach to computing adaptive power flow approximations is scalable to large systems. The power flow calculations for each sample can be performed in parallel and, with good initializations from a nominal operating point, Newton-based power flow solvers will usually converge in a few iterations. These power flow computations can be reused across bus voltage and line flow approximations. The regression problems for calculating the approximations' coefficients are tractable linear or quadratic programs that are independent of each other and can thus also be parallelized. Note also that constraint screening methods can identify redundant limits to reduce the need to compute approximations of some quantities~\cite{owen_aquino_roald_molzahn-ac_constraint_screening}.

We also emphasize that these calculations can be performed offline across a forecasted range of operation and then employed online with linear programming solvers or to make bilevel and mixed-integer problems tractable. Example applications of CLAs include bilevel problems for siting sensors~\cite{buason2023datadriven}, charging optimization for electric vehicles prior to evacuations~\cite{owen_aquino_talkington_molzahn-EVacuation_feeder}, and more fair curtailments of solar photovoltaic generation~\cite{gupta_buason_molzahn-fairness_pv_limit}; see also the survey papers~\cite{10202779, jia2023tutorial1, jia2023tutorial2} on adaptive power flow approximations.}


\section{Numerical results} \label{sec:simulation}
This section validates the proposed importance sampling and rational function approximation methods for the power flow equations. We begin by comparing the [1/1] multivariate generalization of the Pad\'e approximant to the first- and second-order Taylor approximations for voltage magnitudes. Following that, we evaluate the performance of linear approximations (LAs), rational approximations (RAs), as well as their conservative versions, CLAs and CRAs. Additionally, we utilize the second-order sensitivities in the importance sampling. Lastly, we provide a conceptual demonstration through a simplified optimal power flow problem.

The test cases used in the simulations include the IEEE 24-bus system, \textit{case30}, \textit{case33bw}, \textit{case85}, \textit{case141}, and \textit{case2383wp}, all of which are available in M{\sc atpower}~\cite{zimmerman_matpower_2011}. 
For the voltage and current flow approximations, we draw $500-1000$ samples by varying the power injections from $70\%$ to $130\%$ of their nominal values. The voltage and current flow values are reported in per unit (pu). We use $\ell_1$ for a loss function $\mathcal{L}(\,\cdot\,)$. \textcolor{black}{The hardware used for the computations in this section is a MacBook Pro with an Apple M1 Pro chip with 10 cores and 16~GB of RAM.}

\begin{figure}[b!]
\vspace{-1.8em}
	\centering 
    \subfloat[Predicted vs actual voltage\label{1a-taylor_pade}]{\includegraphics[width=0.5\linewidth]{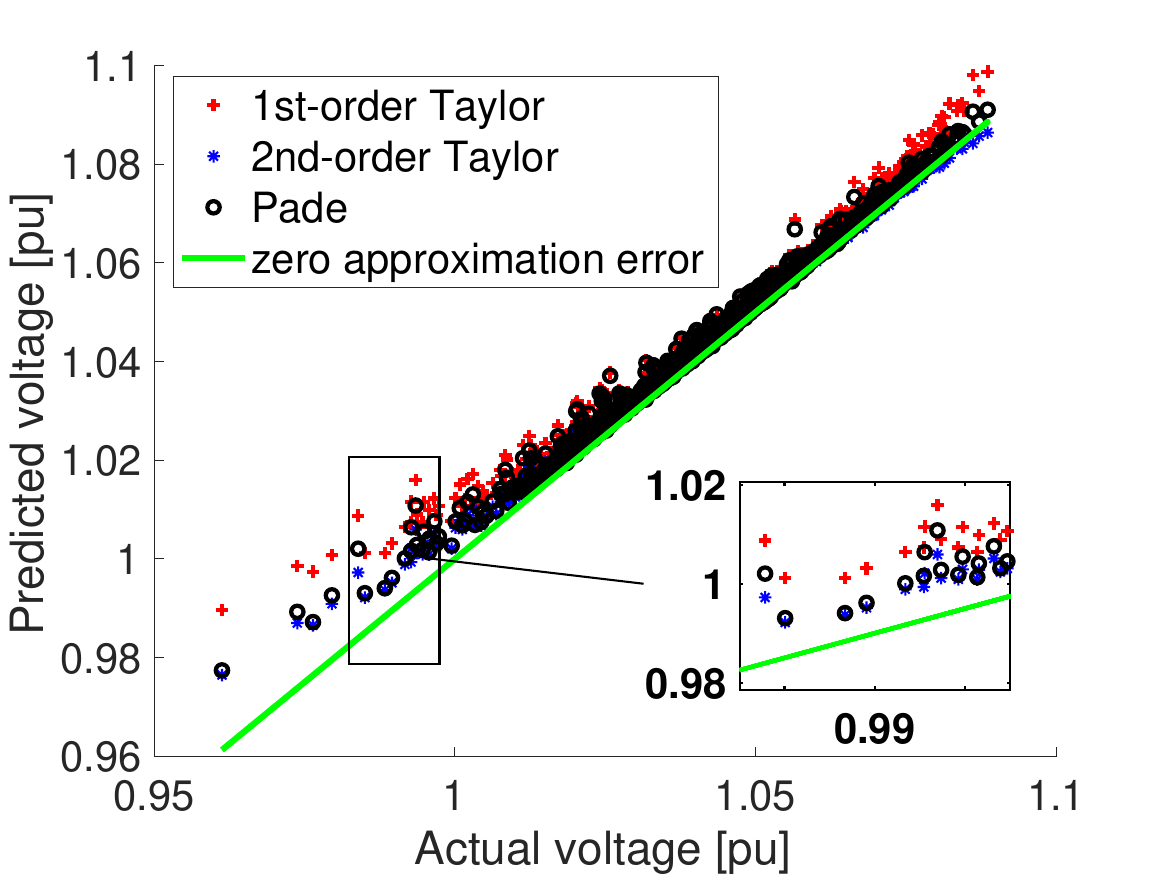}}
    \hfill
    \subfloat[Histogram of errors\label{1b-taylor_pade_histogram}]{\includegraphics[width=0.5\linewidth]{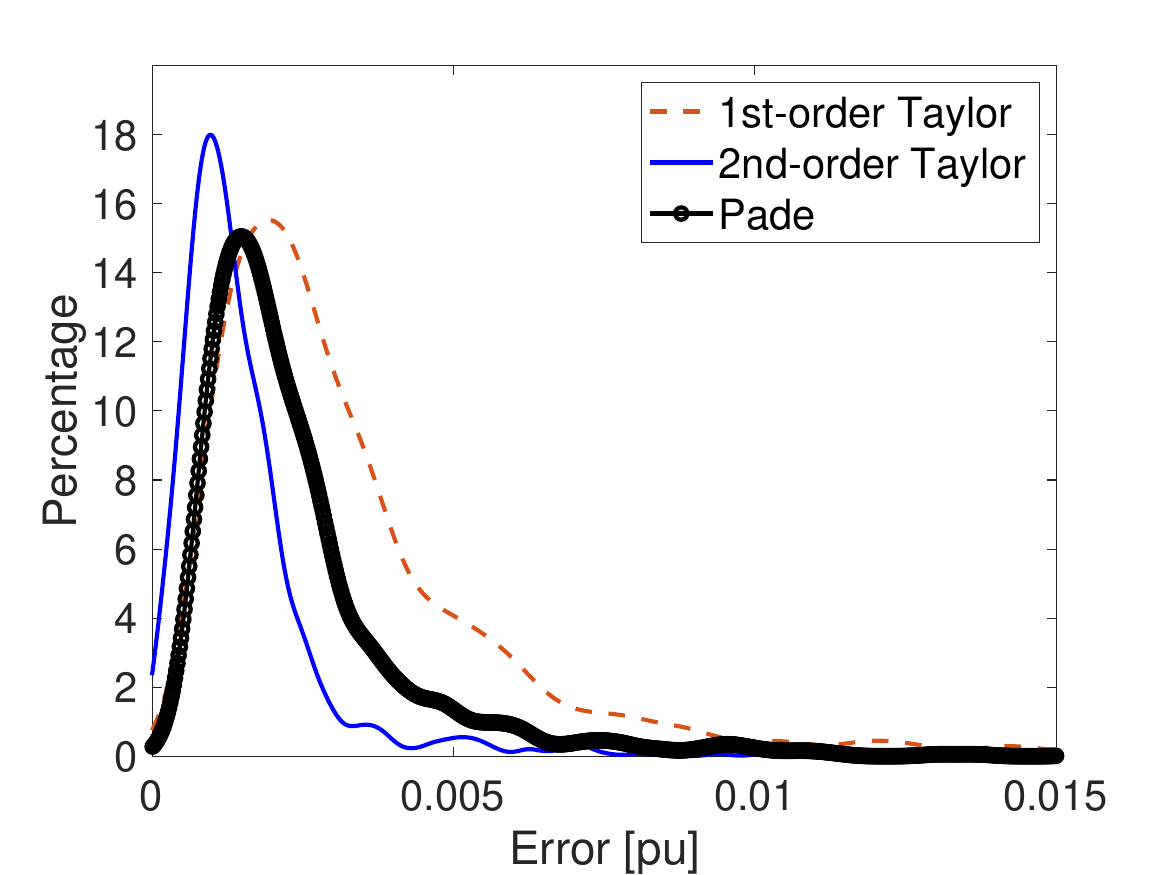}}
	\caption{The left plot shows a comparison between the first- (red crosses) and second-order Taylor approximations (blue asterisks) and the Pad\'e approximant (black circles) for voltage magnitudes at bus 1 in the IEEE 24-bus system. The green line at $45^\circ$ represents zero approximation error. The smoothed histogram in the right plot displays the errors from the first- (red dashed line) and second-order Taylor approximations (blue solid line), along with the [1/1] multivariate generalization of the Pad\'e approximant (thick black line).}
	\label{fig:Pade vs Taylor_ieee 24_bus1}
	\vspace{-0.5em}
\end{figure}


\begin{table}[b!] 
\caption{Voltage magnitude errors at a specific bus} \label{table:pade_taylor}
\centering
\setlength\tabcolsep{2pt}
\resizebox{\columnwidth}{!}{
\begin{tabular}{c|c|c|c|c}
  \multirow{3}{*}{Cases} & \multirow{3}{*}{Bus} & \multicolumn{3}{c}{Average errors per sample [pu]} \\
  \cline{3-5}
  & & $1^{\text{st}}$-order  & [1/1] multivariate & $2^{\text{nd}}$-order \\
  & & Taylor & Pad\'e ($^*$) & Taylor ($^*$) \\
  \hline \hline
  IEEE $24$-bus & 22 & $4.72\times10^{-3}$ & $2.79\times10^{-3}$ (40.9\%) & $9.27\times10^{-4}$ (80.4\%) \\
  \hline
    \textit{case30} & 30  & $1.38\times10^{-3}$ & $9.78\times10^{-4}$ (29.1\%) & $7.6\times10^{-5}$ (94.5\%) \\
  \hline
  \textit{case33bw} & 33  & $3.2\times10^{-5}$ & $1.6\times10^{-5}$ (50\%) & $3.6\times10^{-7}$ (98.9\%) \\
  \hline
  \textit{case141} & 80 & $6.5\times10^{-6}$ & $3.0\times10^{-6}$ (53.8\%) & $3.8\times10^{-8}$ (99.4\%) \\
  \hline
\end{tabular}}
\noindent\footnotesize{*Percentage error reduction compared the first-order Taylor approximation.}
\vspace{-1.5em}
\end{table}

\begin{table*}[t!] 
\caption{Approximation errors for voltage magnitudes at a specific bus}
\label{table:voltage_LA_RA}
\centering
\resizebox{2\columnwidth}{!}{
\setlength\tabcolsep{2.5pt}
\begin{tabular}{c|c|c|c|c|c|c|c|c|c|c|c}
  \multirow{3}{*}{Cases} & \multirow{3}{*}{Bus} & \multicolumn{10}{c}{Errors per sample [pu]} \\
  \cline{3-12}
  & & \multicolumn{2}{c|}{LA} & \multicolumn{2}{c|}{RA} & \multicolumn{3}{c|}{CLA} & \multicolumn{3}{c}{CRA}\\
  \cline{3-12}
  & & Mean & Max & Mean ($^*$) & Max & Mean & Max & time [s] & Mean ($^\dagger$) & Max & time [s] \\
  \hline \hline
    \multirow{2}{*}{\textit{case30}} & \multirow{2}{*}{25} & \multirow{2}{*}{$2.55\times10^{-3}$} & \multirow{2}{*}{$2.76\times10^{-2}$}  & $2.18\times10^{-3}$  & \multirow{2}{*}{$2.14\times10^{-2}$}  & \multirow{2}{*}{$4.81\times10^{-3}$} & \multirow{2}{*}{$3.46\times10^{-2}$} & \multirow{2}{*}{0.465} & $3.33\times10^{-3}$  & \multirow{2}{*}{$1.61\times10^{-2}$} & \multirow{2}{*}{0.608}  \\
    & & & & (14.51\%) & & & & & (30.77\%) & & \\
  \hline
  \multirow{2}{*}{\textit{case33bw}} & \multirow{2}{*}{33} & \multirow{2}{*}{$9.06\times10^{-5}$} & \multirow{2}{*}{$8.40\times10^{-4}$}  & $3.37\times10^{-5}$ & \multirow{2}{*}{$3.81\times10^{-4}$}  & \multirow{2}{*}{$1.37\times10^{-4}$} & \multirow{2}{*}{$8.78\times10^{-4}$} & \multirow{2}{*}{0.471} & $1.21\times10^{-4}$ & \multirow{2}{*}{$9.61\times10^{-4}$} & \multirow{2}{*}{0.611}  \\
  & & & & (62.80\%) & & & & & (11.68\%) & & \\
  \hline
  \multirow{2}{*}{\textit{case141}} &  \multirow{2}{*}{80} &  \multirow{2}{*}{$1.87\times10^{-5}$} &  \multirow{2}{*}{$1.74\times10^{-4}$}  &  $4.01\times10^{-6}$ &  \multirow{2}{*}{$4.94\times10^{-5}$}  &  \multirow{2}{*}{$2.62\times10^{-5}$} &  \multirow{2}{*}{$2.49\times10^{-4}$} &  \multirow{2}{*}{0.635} & $1.91\times10^{-5}$ &  \multirow{2}{*}{$2.09\times10^{-4}$} &  \multirow{2}{*}{0.874} \\
  & & & & (78.56\%) & & & & & (16.06\%) & & \\
  \hline
   \multirow{2}{*}{\textit{case2383wp}} &  \multirow{2}{*}{466} &  \multirow{2}{*}{$5.60\times10^{-6}$} &  \multirow{2}{*}{$2.31\times10^{-5}$} &  $3.79\times10^{-6}$ &  \multirow{2}{*}{$1.57\times10^{-5}$} &  \multirow{2}{*}{$1.02\times10^{-5}$} &  \multirow{2}{*}{$4.52\times10^{-5}$}  &  \multirow{2}{*}{225.2$\textsuperscript{\textsection}$}  & $9.32\times10^{-6}$  &  \multirow{2}{*}{$3.55\times10^{-5}$} &  \multirow{2}{*}{238.3$\textsuperscript{\textsection}$}  \\
   & & & & (32.32\%) & & & & & (8.63\%) & & \\
  \hline
\end{tabular}}
\footnotesize{*The percentage reduction in errors compared to the mean errors from linear approximation (LA). \\ $\dagger$The percentage reduction in errors compared to the mean errors from conservative linear approximation (CLA). \\ $\textsuperscript{\textsection}$Time to solve 1000 power flow solutions is 206.1 s. Times to solve for CLA and CRA are 19.1 and 32.2 s, respectively.}
\end{table*}

\begin{table*}[t!] 
\vspace{-1em}
\caption{Approximation errors for the current flows on a specific line}
\label{table:flow_LA_RA}
\centering
\resizebox{2\columnwidth}{!}{
\setlength\tabcolsep{2.5pt}
\begin{tabular}{c|c|c|c|c|c|c|c|c|c}
  \multirow{3}{*}{Cases} & \multirow{3}{*}{Line} & \multicolumn{8}{c}{Errors per sample [pu]} \\
  \cline{3-10}
  & & \multicolumn{2}{c|}{LA} & \multicolumn{2}{c|}{RA} & \multicolumn{2}{c|}{CLA} & \multicolumn{2}{c}{CRA}\\
  \cline{3-10}
  & & Mean & Max & Mean ($^*$) & Max & Mean & Max & Mean ($^\dagger$) & Max \\
  \hline \hline
    \textit{case30} & 1-2  & $1.25 \times 10^{-2}$ & $8.82 \times 10^{-2}$ & $1.20 \times 10^{-2}$ (4\%) & $8.56 \times 10^{-2}$ & $4.28 \times 10^{-2}$ & $1.48 \times 10^{-1}$ & $3.63 \times 10^{-2}$ (15.19\%) & $1.34 \times 10^{-1}$  \\
  \hline
  \textit{case33bw} & 29-30  & $2.61 \times 10^{-4}$ & $2.15 \times 10^{-3}$ & $2.49 \times 10^{-4}$ (4.60\%) & $1.69 \times 10^{-3}$ & $9.32 \times 10^{-4}$ & $3.25 \times 10^{-3}$ & $5.68 \times 10^{-4}$ (39.06\%) & $1.95 \times 10^{-3}$ \\
  \hline
    \textit{case85} & 3-17 & $9.23 \times 10^{-4}$ & $7.41 \times 10^{-3}$ & $8.33 \times 10^{-4}$ (9.75\%) & $6.51 \times 10^{-3}$ & $2.81 \times 10^{-3}$ & $9.10 \times 10^{-3}$ & $1.92 \times 10^{-3}$ (31.67\%) & $8.17 \times 10^{-3}$ \\
  \hline
  \textit{case141} & 92-93 & $3.03 \times 10^{-4}$ & $3.76 \times 10^{-3}$ & $1.73 \times 10^{-4}$ (42.90\%) & $2.01 \times 10^{-3}$ & $6.98 \times 10^{-4}$ & $3.09 \times 10^{-3}$  & $2.57 \times 10^{-4}$ (63.18\%) & $3.42 \times 10^{-3}$ \\
  \hline
\end{tabular}}
\footnotesize{*The percentage reduction in errors compared to the mean errors from linear approximation (LA). \\ $\dagger$The percentage reduction in errors compared to the mean errors from conservative linear approximation (CLA).}
\end{table*}

\subsection{Multivariate generalization of the Pad\'e approximant} \label{sub:simulation_Pade}

By incorporating curvature information from the second-order Taylor approximation, the [1/1] multivariate generalization of the Pad\'e approximant can surpass the accuracy of the first-order Taylor approximation. We empirically demonstrate these accuracy advantages by comparison to both the first- and second-order Taylor approximation across a range of operation around a nominal point. For the sake of illustration, we choose, in each test case, the voltage magnitude at a bus for which the first-order Taylor approximation has a large error.

Fig.~\ref{fig:Pade vs Taylor_ieee 24_bus1} presents an illustrative example of voltage approximation at bus 1 in the IEEE 24-bus system. The results depicted in Fig.~\ref{fig:Pade vs Taylor_ieee 24_bus1}(a) indicate that the [1/1] multivariate generalization of the Pad\'e approximant has approximation accuracy between the first- and second-order Taylor approximations. To provide a clearer view of the errors associated with each approximation, we present a histogram plot in Fig.~\ref{fig:Pade vs Taylor_ieee 24_bus1}(b).

To quantify these results, Table~\ref{table:pade_taylor} presents the average errors per sample in voltage magnitudes resulting from three different approximation methods: the first-order Taylor approximation, the [1/1] multivariate generalization of the Padé approximant, and the second-order Taylor approximation. The findings reveal substantial error reductions, ranging from 29.1\% to 53.8\%, when employing the [1/1] multivariate Padé approximant in comparison to the first-order Taylor approximation. The second-order Taylor approximation yields larger error reductions, ranging from 80.4\% to 99.4\%, compared to the first-order Taylor approximation, but unlike the [1/1] multivariate Padé approximant, it does not maintain linearity (or even convexity) of constraints when used for optimization. \textcolor{black}{Notably, while the [1/1] Pad\'e approximant matches the second-order Taylor approximation's curvature in the univariate case, this equivalence does not hold in the multivariate version, as demonstrated by the numerical tests.}


\subsection{Linear and rational approximations} \label{sub:simulation_cra}
Next, we compare the performance of four approximation methods: 1) linear approximation (LA), 2) rational approximation (RA), 3) conservative linear approximation (CLA), and 4)~conservative rational approximation (CRA). RA and CRA are rational functions with linear numerators and denominators. The LAs and RAs are similar to the CLAs and CRAs, respectively, except that they do not enforce conservativeness.

\begin{figure}[t]
\vspace{-2em}
	\centering 
	\subfloat[Predicted vs actual flow\label{1a-la_ra}]{\includegraphics[width=0.5\linewidth]{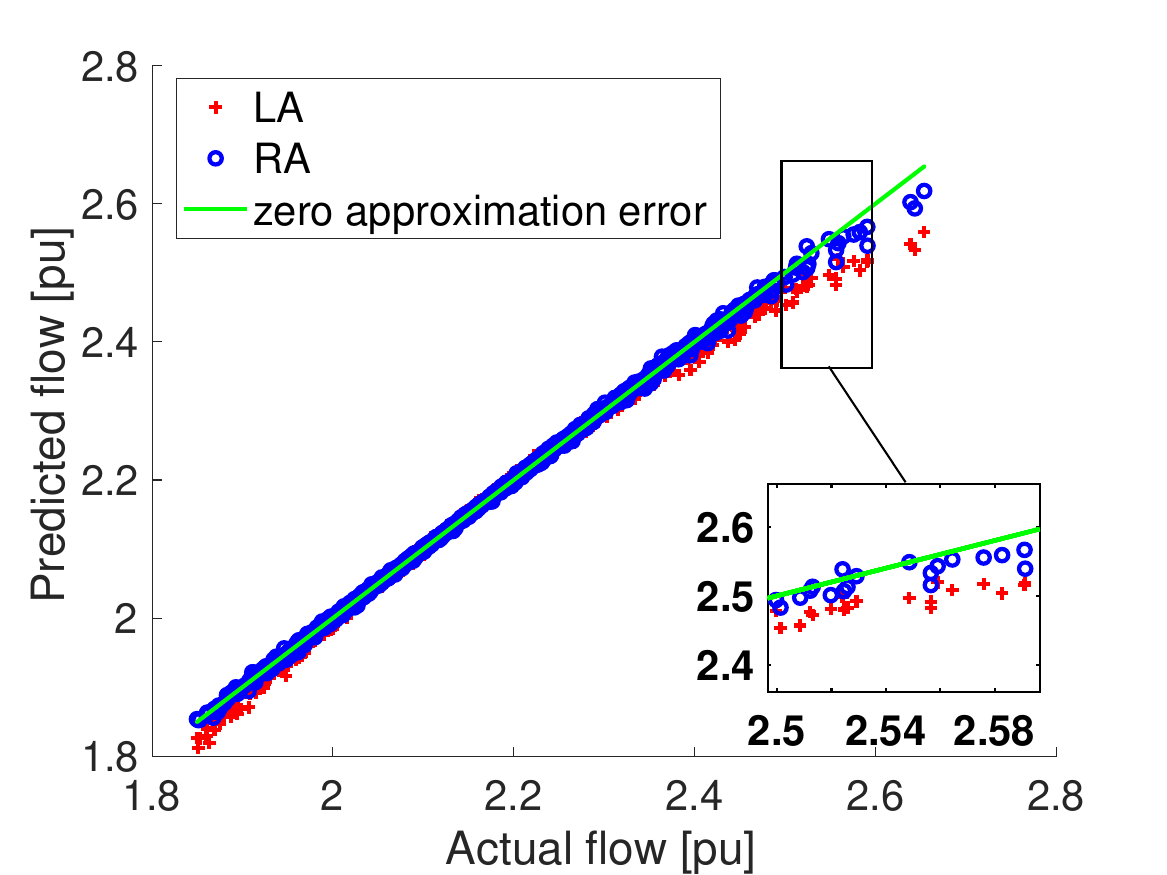}}
    \hfill
  \subfloat[Histogram of errors\label{1b-la_ra_histogram}]{\includegraphics[width=0.5\linewidth]{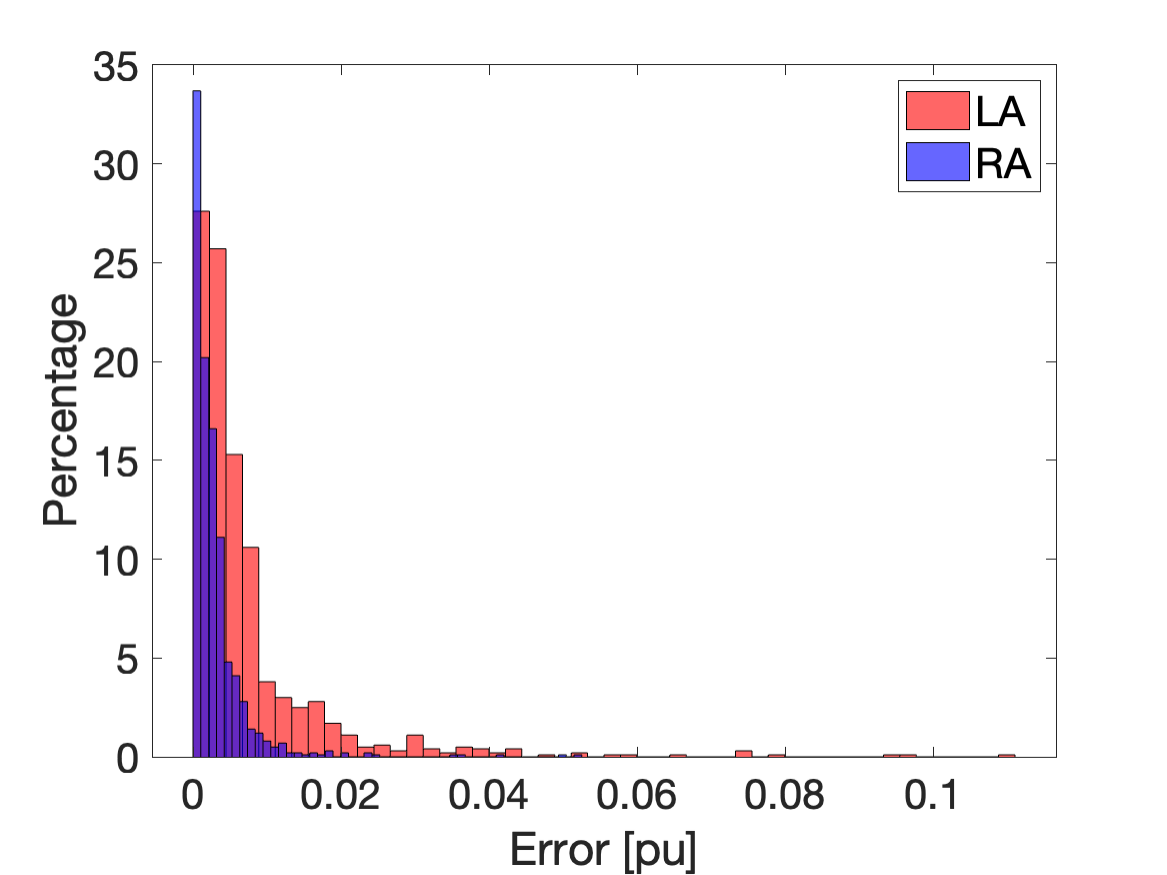}}
	\caption{The left plot compares linear approximation (LA) in red crosses and rational approximation (RA) in blue circles for current flow on a branch connecting buses 15 and bus 21 in IEEE 24-bus system. The green line indicates the zero approximation error. The right plot shows the error histogram from the linear approximation (LA) and the rational approximation~(RA).}
	\label{fig:la vs ra_ieee 24}
	\vspace{-1em}
\end{figure}

In Fig.~\ref{fig:la vs ra_ieee 24}, we provide an illustrative example comparing the performance of LA and RA for current flow from bus 15 to bus 21 in the IEEE 24-bus system. Fig.~\ref{fig:la vs ra_ieee 24}(a) clearly demonstrates that RA predicts the current flow much more accurately in comparison to LA. For a more detailed inspection of the errors, Fig.~\ref{fig:la vs ra_ieee 24}(b) presents a histogram plot. The RA yields significantly lower errors as evident from the left side in the histogram, with very few samples exhibiting errors in excess of $0.01$ pu.

\textcolor{black}{Table~\ref{table:voltage_LA_RA} compares the performance of RA and CRA with LA and CLA, respectively, in terms of approximating voltage magnitudes. In this analysis, load values vary from 30\% to 170\% of their nominal values. The results indicate significant error reductions for RAs, with improvements ranging from 14.51\% to 78.56\% compared to LAs. Similarly, CRAs exhibit error reductions over CLAs, with improvements ranging from 8.63\% to 30.77\%. Regarding computation time, it takes less than a second to compute 500 power flow solutions and to solve for CLAs and CRAs for a specific bus in the \textit{case30}, \textit{case33bw}, and \textit{case141} systems. For the larger \textit{case2383wp} system, computing 1000 power flows takes 206.1 seconds, while solving for CLA and CRA coefficients for a particular bus voltage takes 19.1 and 32.2 seconds, respectively.}\footnote{\textcolor{black}{With 2383 buses, there are a total of 4766 voltage approximations that could be computed for this system. Since each can be computed in parallel, the total computation time depends on the number of processors available. For instance, a 100-processor cluster where each computing node had similar capabilities as the laptop used in these numerical tests could compute all voltage magnitude approximations within approximately fifteen minutes.}} 

Table~\ref{table:flow_LA_RA} presents a similar comparison for current flow approximations. RAs provide error reductions ranging from 4\% to 42.90\% compared to LAs. CRAs achieve even more significant error reductions, ranging from 15.19\% to 63.18\% improvements compared to CLAs. These results demonstrate the advantages inherent to rational approximations.


\subsection{Second-order sensitivity: Computation time} \label{sub:computation time_CLA}
This section presents results on the computation time required for calculating the second-order sensitivity information. 
The computational process is divided into three distinct steps: $a)$ power flow solution, $b)$ second-order sensitivity matrix calculation, and $c)$ singular value decomposition (SVD). 

\begin{table}[ht!]
\caption{Computation time for the second-order sensitivity matrix}
\label{table:second_order_time}
\centering
\begin{tabular}{c|c|c|c}
  & \multicolumn{3}{c}{Computation time [\SI{}{\second}]} \\
  \cline{2-4}
  & Step $a)$ & Step $b)$ & Step $c)^*$ \\
  \hline \hline
  IEEE $24$-bus & 0.0054 & \phantom{1}0.021\phantom{0} & 0.0002\\
  \hline
    \textit{case30} & 0.0040 & \phantom{1}0.022\phantom{0} & 0.002\phantom{0} \\
  \hline
  \textit{case33bw} & 0.0042 & \phantom{1}0.026\phantom{0} & 0.001\phantom{0} \\
  \hline
  \textit{case141} & 0.0029 & \phantom{1}0.1185 & 0.0036  \\
  \hline
  \textit{case2383wp} & 0.048\phantom{0} & 13.07\phantom{00} & 3.58\phantom{00} \\
  \hline
\end{tabular}\newline
\noindent \hspace*{-3em} \footnotesize{*Computation time for finding all singular values.}
\end{table}

The computation times presented in Table~\ref{table:second_order_time} are specific to the analysis of a single bus for each test case. However, note that Step~$a)$ and most parts of Step~$b)$ (e.g., computing Jacobian and Hessian matrices) only need to be computed once for all output quantities, which can result in time savings and increased efficiency during the overall process if multiple output quantities are of interest (e.g., voltage magnitudes at several buses). Additionally, each case shows only a few significant singular values (at least 10\% of the maximum singular value), e.g., three for \textit{case33bw} at bus 18 and three for \textit{case2383wp} at bus 466. Consequently, we can speed up the computation in Step $c)$ by evaluating only a limited number of the most significant singular values.

\subsection{Second-order sensitivity: Span of singular vectors} \label{sub:sim_span}
The total time for computing sample-based power flow approximations depends on the number of sampled sets of power injections. Using second-order sensitivities, the importance sampling approach described in Section~\ref{sec:adaptive sampling} aims to reduce computing time by prioritizing power injection samples in regions of high curvature. If there were substantial variability in the second-order sensitivities across typical operating ranges, this importance sampling approach would likely perform poorly as the second-order sensitivity information at a subset of samples would not help characterize the entire operating range. We next provide empirical evidence showing that the space spanned by the dominant singular values change very little over wide operating ranges.

We consider 1000 randomly sampled power injections that range from 70\% to 130\% of their nominal values and compute the first $k$ significant singular values ($k=3$ for \textit{case33bw} and $k=3$ for \textit{case2383wp}) for the second-order sensitivity matrix for each of these points. We stack the corresponding singular vectors to form a matrix, denoted as $\mathbf{M}$, consisting of $1000k$ columns and compute the singular values of this matrix. The number of significant singular values of $\mathbf{M}$, i.e., its approximate rank, indicates the extent to which the space spanned by the dominant singular vectors vary across the power flow samples. If $M$ has approximate rank $k$, this indicates that the  singular vectors have no variation, whereas if $\mathbf{M}$ has full rank then the singular vectors vary significantly over the operating range.

\begin{figure}[ht!] 
\vspace{-1.75em}
    \centering
  \subfloat[\textit{case33bw} at bus 18\label{1a-svd_1}]{\includegraphics[width=0.5\linewidth]{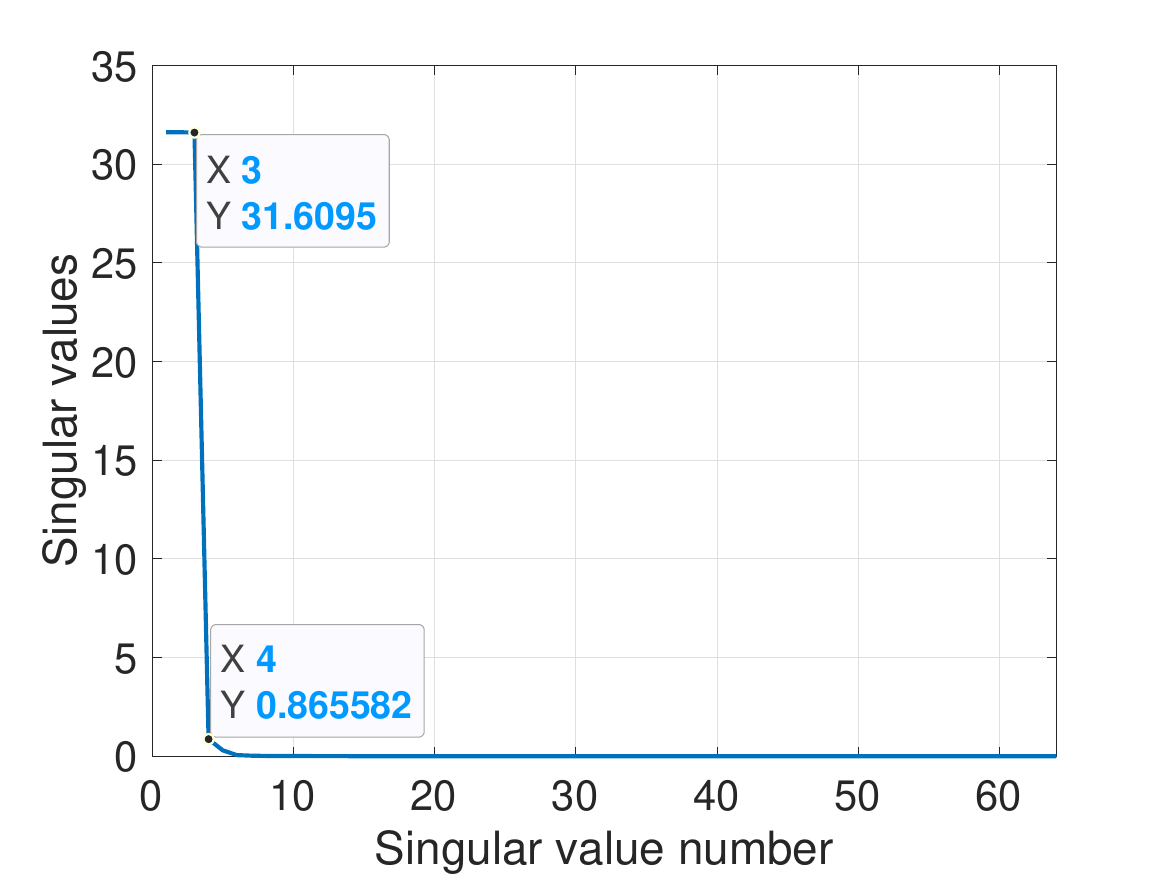}}
    \hfill
  \subfloat[\textit{case2383wp} at bus 466\label{1b-svd_2}]{\includegraphics[width=0.5\linewidth]{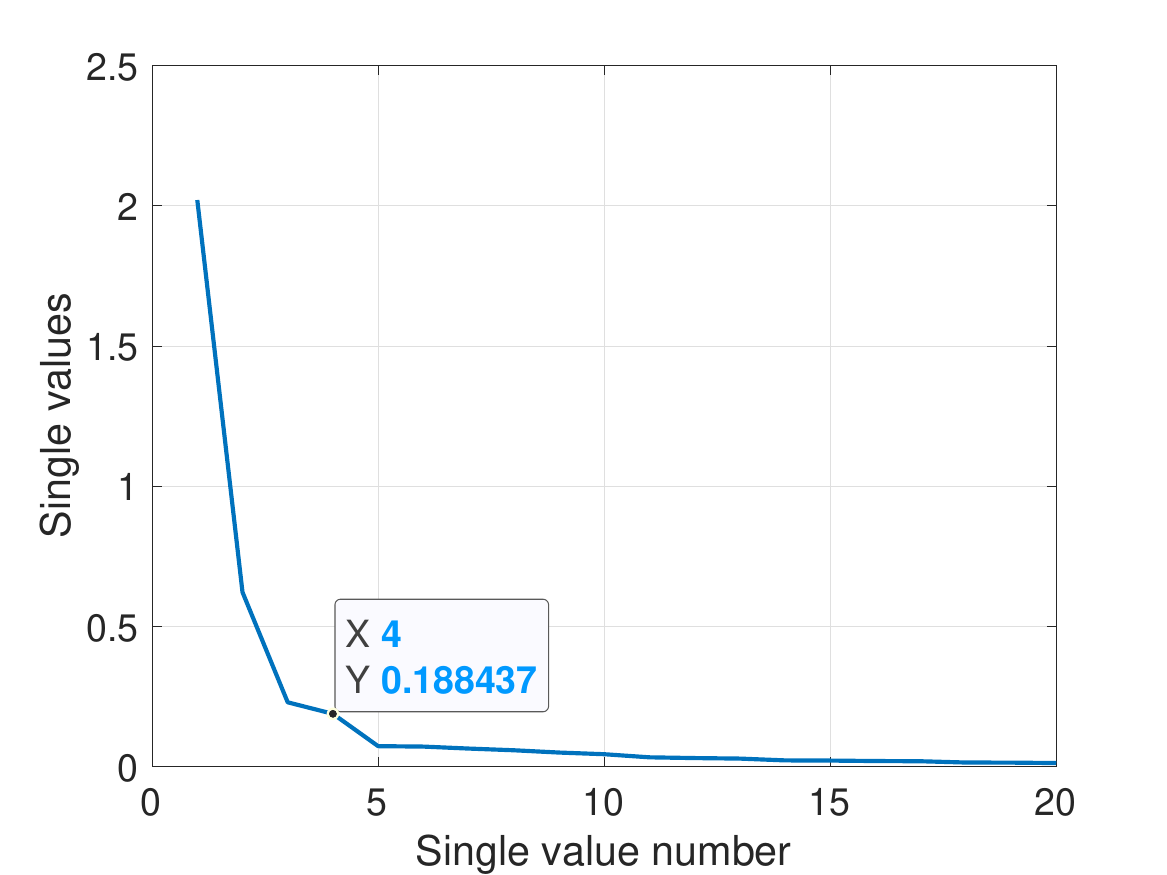}}
   \caption{The singular values of $\mathbf{M}$ for 1000 spans of the singular vectors. Fig.~(b) only shows first 20 singular values.}
  \label{fig:singular_values} 
  \vspace{-1.5em}
\end{figure}

As illustrative examples, Fig.~\ref{fig:singular_values} plots the sorted singular values of $M$ for \textit{case33bw} and \textit{case2383wp}. The plot reveals that for the \textit{case33bw} ($k=3$), the stacked matrix $\mathbf{M}$ has three large singular values, followed by two non-negligible but considerably smaller ones. The singular values beyond the eighth position for both test cases are negligible. The same behaviour is observed for \textit{case2383wp} with $k=3$. From the results on these and other test cases, we conclude that the second-order sensitivities at a single nominal point is sufficient to characterize the second-order behaviour of the power flow manifold for large ranges of power injections.

\subsection{Second-order sensitivity: Concavity} \label{sub:eigenvalues}
To further analyze the curvature of the power flow manifold, we examine the eigenvalues of the second-order sensitivity matrix. Since the sensitivity matrix is symmetric, all eigenvalues are real-valued. If all eigenvalues are positive, the second-order sensitivity matrix is positive definite, indicating local convexity of the second-order approximation. Conversely, all eigenvalues being negative would indicate local concavity of the second-order approximation. Otherwise, the second-order approximation is indefinite.

\begin{table}[b!]
\vspace{-1.5em}
\caption{Maximum and minimum eigenvalues of the second-order sensitivity matrices for voltage magnitudes}
\label{table:eigenvalues}
\centering
\begin{tabular}{c|c|c|c}
  \multirow{2}{*}{Cases} & \multirow{2}{*}{Bus} & \multicolumn{2}{c}{Eigenvalues} \\
  \cline{3-4}
  & & Max & Min  \\
  \hline \hline
  IEEE $24$-bus  &7&  $-0.002$ & $\phantom{1}{-}0.72\phantom{0}$  \\
  \hline
    \textit{case30} & 30  & $-0.009$ & $\phantom{1}{-}3.02\phantom{0}$  \\
  \hline
  \textit{case33bw} & 18 & $0$ & ${-}10.45\phantom{0}$  \\
  \hline
  \textit{case85} & 50 & $0$ & $\phantom{1}{-}0.34\phantom{0}$  \\
  \hline
  \textit{case141} & 52  & $0$ & $\phantom{1}{-}0.65\phantom{0}$  \\
  \hline
     \textit{case2383wp} & 466 & $0$ & $\phantom{1}{-}2.44\phantom{0}$  \\
  \hline
\end{tabular}
\vspace{-1.5em}
\end{table}

Table~\ref{table:eigenvalues} presents the maximum and minimum eigenvalues of the second-order sensitivity matrices. In each case, we select the bus characterized by the most extreme curvature, which corresponds to the smallest eigenvalue. As suggested by the results in Section~\ref{sub:sim_span}, the second-order sensitivity matrices at a bus exhibit similarity across a range of operational conditions. Therefore, the eigenvalues for each test case are computed at a nominal value. These results empirically demonstrate that all eigenvalues of the second-order sensitivity matrix are non-positive, suggesting local concavity of the second-order approximation. Based on this concavity, we bias the importance sampling procedure to draw more samples at extreme points (i.e., away from the nominal point) for underestimating approximations and around the nominal points for overestimating approximations. The implications of these results will be further discussed in the next section. 

\subsection{Importance sampling} \label{sub:sim_violation}

\begin{figure*}[t!]
	\centering
  \subfloat[All uniform random samples\label{1a-svd}]{\includegraphics[width=0.30\linewidth]{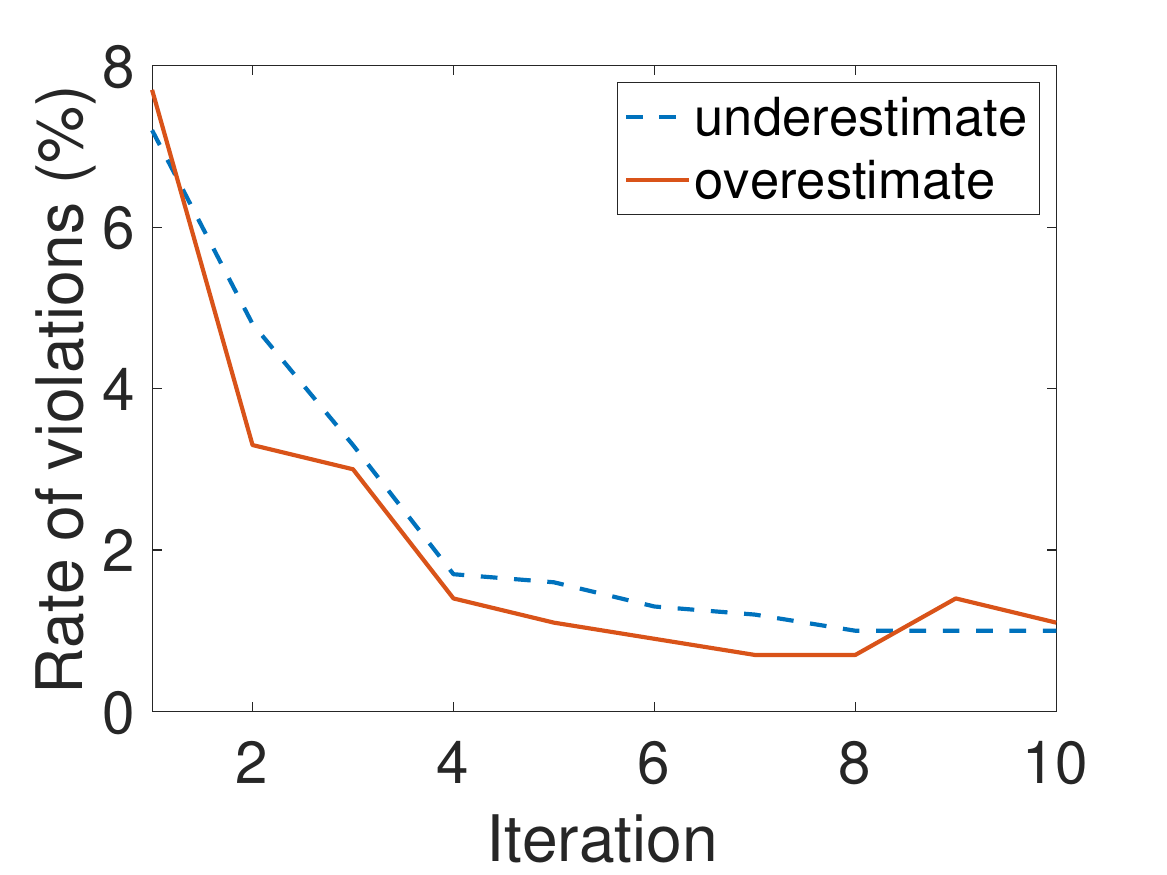}}
    \hfill
  \subfloat[All samples in the dominant singular vectors' span\label{1b-svd}]{\includegraphics[width=0.30\linewidth]{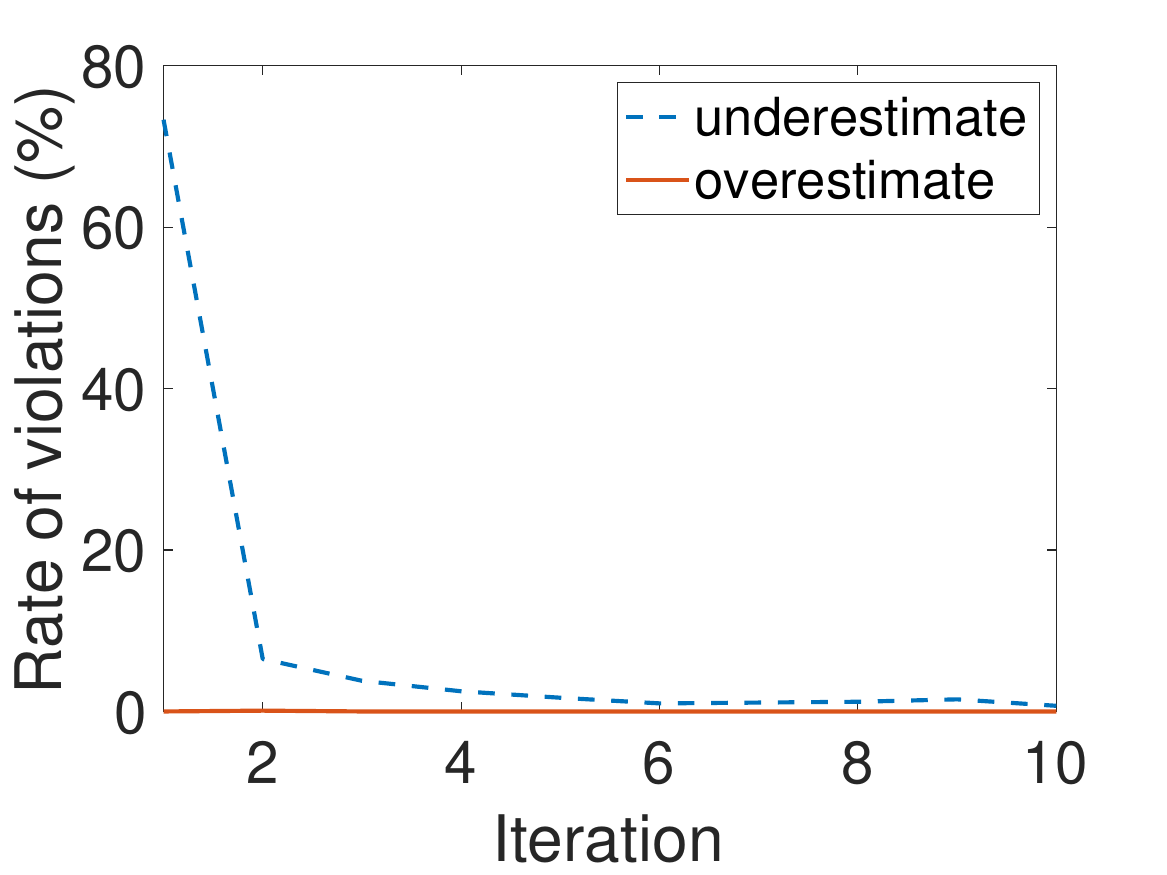}}
  \hfill
  \subfloat[Half uniform random samples and half samples in the dominant singular vectors' span\label{1c-svd}]{\includegraphics[width=0.30\linewidth]{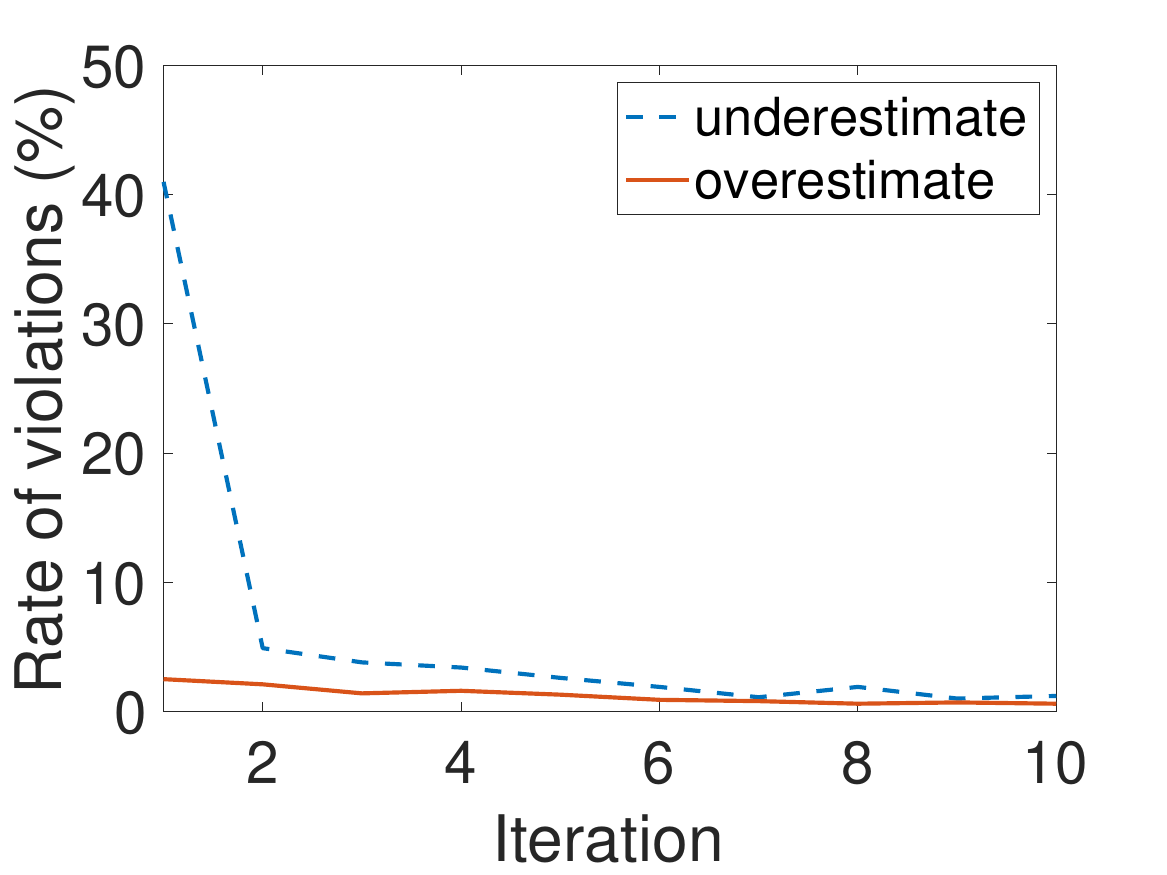}}
   \caption{Plots depicting the rate of violations versus iterations obtained by employing three different sampling strategies: (a) drawing all additional samples randomly, (b) selecting all additional samples in alignment with the subspace spanned by the first three singular vectors, and (c) one half of the additional samples randomly and the other half in accordance with the subspace spanned by the first three singular vectors. These simulations are the approximations of the voltage magnitude at bus $33$ for \textit{case33bw}. The dashed blue and solid red lines denote the violation rates for underestimates and overestimates, respectively.}
  \label{fig:violation_case33} 
  \vspace{-1.5em}
\end{figure*}

Our prior work in~\cite{BUASON2022} introduced a sample selection method aimed at enhancing the conservativeness of the conservative linear approximations with limited impacts on the computation time for solving the regression problem~\eqref{eq:regression}. The sample selection approach in~\cite{BUASON2022} is based on new samples drawn uniformly at random within the predefined range of power injections, which may not adequately capture the non-linearity of the power flow equations or effectively address high curvature regions where more violations can occur.

\begin{figure}[t] 
\vspace{-0.3em}
    \centering
  \subfloat[Underestimate\label{1a-ss}]{\includegraphics[width=0.5\linewidth]{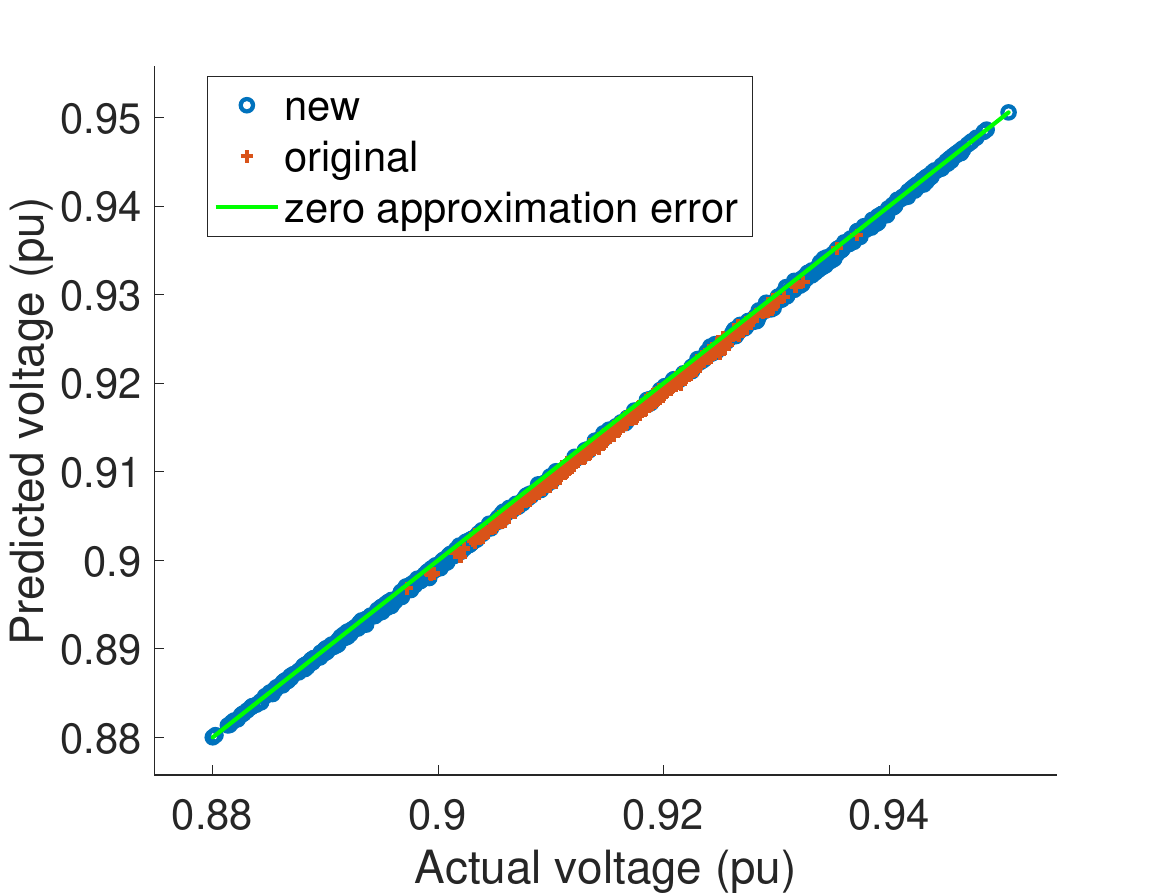}}
    \hfill
  \subfloat[Overestimate\label{1b-ss}]{\includegraphics[width=0.5\linewidth]{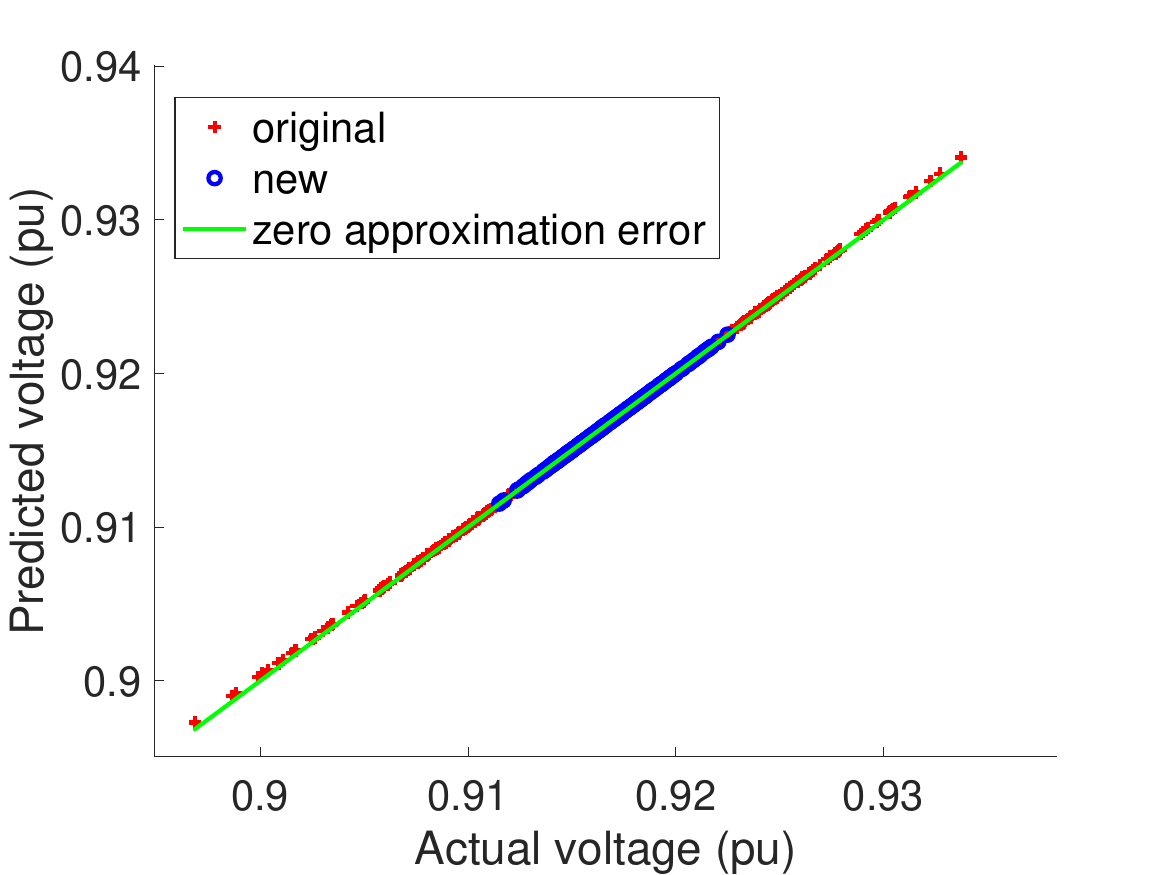}}
   \caption{Plots of the results from (a) underestimating and (b) overestimating linear approximations of the voltage magnitude at bus $33$ for \textit{case33bw}. The red and blue points represent the original and new samples, respectively. The new samples are those that violate the original linear approximations. The green line at $45^\circ$ represents zero approximation error.}
  \label{fig:predicted_voltage_svd_case33} 
  \vspace{-1.5em}
\end{figure}

To address these limitations, our proposed importance sampling approach selects new samples along the span of the dominant singular vectors of the second-order sensitivity matrix. As in~\cite{BUASON2022}, we iterate between selecting new samples of power injections and updating the linearization based on the samples that are not conservative (are above the overestimating function or below the underestimating function). However, as discussed in Section~\ref{sec:adaptive sampling}, we also select samples along the span of the dominant singular vectors to prioritize directions associated with the largest curvature of the power flow equations. 

We also examine the specific neighborhoods from which we draw samples. This later approach takes into account the insights gained from the convexity/concavity analysis in Section~\ref{sub:eigenvalues}. By considering the influence of these directions, we aim to exploit the power flow equations' nonlinear behavior as well as their local convexity or local concavity. 

To assess this importance sampling approach, Fig.~\ref{fig:violation_case33} shows the rates of violations where the newly drawn samples' predicted values are below the approximation's underestimates or exceed the overestimates. When new samples are drawn uniformly at random (Fig.~\ref{fig:violation_case33}(a)), the first iteration's violation rates for both underestimates and overestimates range from 7\% to 8\%, meaning that relatively little new information is provided by additional samples with respect to the conservativeness of the approximation. However, when new samples are drawn in the direction of the subspace spanned by the three singular vectors corresponding to the three largest singular values (Fig.~\ref{fig:violation_case33}(b)), there is a 72\% violation rate for underestimates, but the violation rate for overestimates remains close to 0\% in the first iteration. Violation rates significantly reduce after the first iteration and continue to decrease over the next iterations. These results suggest that, as expected, samples in the span of the dominant singular vectors effectively identify violations for underestimating approximations but have little value for the conservativeness of overestimating approximations. 

To balance the quality of the over- and under-estimating approximations, we assess performance when one half of the additional samples are drawn uniformly at random and the the other half is selected in the direction of the dominant singular vectors. As shown in Fig.~\ref{fig:violation_case33}(c), the violation rate is 41\% for underestimates and 3\% for overestimates, which suggests that a mix of sampling approaches using the second-order sensitivity matrix can substantially improve the performance of the underestimating approximations without overly detrimental impacts to overestimating approximations.

Finally, to provide insights into the characteristics of the violated samples, Fig.~\ref{fig:predicted_voltage_svd_case33} shows both the under- and over-estimating linear approximations of voltages after the importance sampling method (half uniform random samples and the rest in the dominant singular vectors’ span). Consistent with the eigenvalues, the voltage magnitude, when considered as a function of power injections, has local concavity in its second-order curvature. As shown in Figure~\ref{fig:predicted_voltage_svd_case33}(a), which presents the underestimating linear approximation, the newly drawn samples that violate the original CLA are mainly positioned away from the middle of the voltage range. Conversely, the violated samples are concentrated in middle of the voltage range for the overestimating case (see Fig.~\ref{fig:predicted_voltage_svd_case33}(b)). These results emphasize the concavity characteristics discussed in Section~\ref{sub:eigenvalues}.

\pap{In summary, the main goals of the importance sampling method are to capture the non-linearity in power flow equations and focus on high-curvature regions where violations are likely. Our results in Fig.~\ref{fig:violation_case33} serves as an indicator of the nonlinearity in the power flow equations. Uniform random sampling may miss some high-curvature directions, leading to potential violations. The violation rate converges more quickly when using importance sampling compared to random sampling. By using second-order sensitivities, we better sample along nonlinear directions. Additionally, Fig.~\ref{fig:predicted_voltage_svd_case33} shows that the overestimating CLA benefits from samples near the middle of the operating range, while the underestimating CLA gains from samples farther away (also illustrated in Fig.~\ref{fig:overestimate_plane_example}).}

\subsection{Application: Simplified optimal power flow} \label{sub:sim_opf}

While the primary intended applications of our RA and CRA methods are bilevel problems~\cite{wogrin2020,buason2023datadriven} and capacity expansion planning problems~\cite{Go2016}, demonstration on these problems is beyond the scope of this paper. Instead, we focus on presenting results in a simplified optimal power flow (OPF) setting, which serves as a conceptual demonstration and a basis for comparison between various linear approximations (DC power flow, LA, CLA, RA, and CRA). The simplified versions of the OPF problem enforce constraints on voltages at load buses (where $P$ and $Q$ are specified) and power generations within specified limits, without considering line flow limits.

For LA, CLA, RA, and CRA, we compute approximations for voltages at load buses, reactive power outputs at generator buses (where $P$ and $V$ are specified), and active power generation at a reference bus. These approximations are functions of active power injections and voltages at generator buses.

\begin{table}[t]
\caption{Results comparing solutions from\\ AC-, DC-, LA-, CLA-, RA-, and CRA-OPF}
\label{table:opf}
\vspace{-0.5em}
\centering
\begin{tabular}{ c|c|c|c }
& \textit{case6ww} 
& \textit{case9}   & \textit{case14} 
 \\ \cline{2-4} 
\hline \hline
\textbf{AC-OPF}       & 2986.04  & 1456.83 & 5368.30       \\ \hline
\textbf{DC-OPF}     & 2995.15 (0.31\%)    & 1502.82 (3.16\%)  & 5368.52 (0.004\%)\\ 
   \textit{Violation} & $V$ (0.029 pu) & -  & - \\
   \hline
\textbf{LA-OPF}   & 2989.95 (0.13\%)   & 1473.67 (1.16\%)  & 5368.52 (0.004\%)\\
 \textit{Violation} & - & - & $V$ (0.004 pu) \\
\hline             
\textbf{CLA-OPF}   & 2991.14 (0.17\%)  & 1475.92 (1.31\%) & 5368.52 (0.004\%)\\
\textit{Violation} & - & -  & - \\ 
\hline
\textbf{RA-OPF}   & 2988.93 (0.10\%)    & 1471.25 (0.99\%)  &  5368.51 (0.004\%)\\
 \textit{Violation} & - & - & $V$ (0.002 pu) \\
\hline             
\textbf{CRA-OPF}   & 2990.25 (0.14\%)  & 1473.91 (1.17\%)  & 5368.51 (0.004\%)\\
\textit{Violation} & - & -  & - \\ 
\hline
\end{tabular}
\footnotesize{The values in parentheses $(\,\cdot\,)$ indicate the percentage difference in cost when compared to the cost obtained from AC-OPF.}
\vspace{-2em}
\end{table}

Table~\ref{table:opf} presents a comparison of results obtained from various approximations with AC-OPF for \textit{case6ww}, \textit{case9}, and \textit{case14}. In each cell of this table, the first row displays the actual cost associated with the AC-PF feasible solution, determined by the set points prescribed in different OPF formulations. The second row indicates the presence of voltage violations and quantifies the maximum voltage violation.

The results show that both CLA-OPF and CRA-OPF do not result in any voltage violations due to their conservativeness property. In \textit{case6ww}, DC-OPF leads to a maximum voltage violation of 0.029 pu. For \textit{case14}, LA-OPF and RA-OPF produce a maximum voltage violation of 0.002 pu and 0.004 pu, respectively. Moreover, RA-OPF exhibits a cost advantage over LA-OPF and DC-OPF in \textit{case6ww} and \textit{case9}. Similarly, CRA-OPF demonstrates cost improvements in comparison to CLA-OPF and DC-OPF in these cases.

\section{Conclusion and future work} \label{sec:future work}
\par This paper presents an importance sampling approach designed to enhance the accuracy and conservativeness (i.e., the tendency to over- or under-estimate a quantity of interest) of power flow approximations. Our method leverages second-order sensitivity information to provide a deeper understanding of the relationships between various quantities. Additionally, inspired by the Pad\'e approximant and second-order sensitivities, we introduce the [1/1] multivariate Pad\'e approximant, expressed as a rational function with a linear numerator and denominator, enhancing accuracy beyond the capabilities of linear functions. This rational approximation, when used as a constraint, preserves linearity in decision variables. Our numerical results reveal the benefits of second-order sensitivities for the importance sampling method and demonstrate improved accuracy compared to other linear approximations.

Our future work aims to develop second-order sensitivities for various output functions, such as $V^2$, and current flow. We will also explore the [2/2] multivariate Pad\'e approximant to construct convex quadratic approximations. Additionally, we will focus on the applications of our proposed method to power system planning and resilience tasks, including capacity expansion planning problems. \textcolor{black}{We also plan to develop linear and rational approximations applicable to multiple topologies, where certain lines may be switched on or off.}

\section*{Appendix: Second-Order Sensitivity Analysis} \label{appendix:sos}

The appendix formulates the second-order sensitivities of the voltage magnitudes with respect to the active and reactive power injections.
Recall the length-$2N$ vectors $x$ (of active and reactive power injections) and $y$ (of voltage angles and magnitudes) from~\eqref{eq:x_y_def}.
We first compute the first-order sensitivities by evaluating the gradient of each row of the vectors in~\eqref{eq:pf_simplified} with respect to $\bm{x}$. The derivative expression is:
\begin{subequations} \label{eq:first derivative}
    \begin{align}
        \bm{I} &= (\bm{J})(\nabla_{\bm{x}} \bm{y}), \label{eq:first_sensitive} \\ 
        \nabla_{\bm{x}} \bm{y} &:= \begin{bmatrix}
            \cfrac{\partial y_1}{\partial x_1} &\cfrac{\partial y_1}{\partial x_2} &\cdots & \cfrac{\partial y_1}{\partial x_{2N}} \\
            \vdots &\vdots &\ddots &\vdots\\
             \cfrac{\partial y_{2N}}{\partial x_1} &\cfrac{\partial y_{2N}}{\partial x_2} &\cdots & \cfrac{\partial y_{2N}}{\partial x_{2N}}
        \end{bmatrix}, \label{eq:gradient}
    \end{align}
\end{subequations}
where $\bm{I}$ represents an identity matrix with appropriate dimensions and $\bm{J}$ is the Jacobian matrix associated with the power flow equations, expressed as:
\begin{equation}
    \bm{J} = \begin{bmatrix}
            \cfrac{\partial x_1}{\partial y_1} &\cfrac{\partial x_1}{\partial y_2} &\cdots & \cfrac{\partial x_1}{\partial y_{2N}} \\
            \vdots &\vdots &\ddots &\vdots\\
             \cfrac{\partial x_{2N}}{\partial y_1} &\cfrac{\partial x_{2N}}{\partial y_2} &\cdots & \cfrac{\partial x_{2N}}{\partial y_{2N}}
     \end{bmatrix}. \label{eq:Jacobian}
\end{equation}

From~\eqref{eq:first derivative}, $\frac{\partial y_k}{\partial x_i}$ is simply $[\bm{J}^{-1}]_{ki}$, i.e., the $(k,i)$ entry of the inverse Jacobian matrix. The second-order sensitivities are computed by evaluating the partial derivative of each element in the expression for the first-order sensitivities~\eqref{eq:gradient}. The partial derivatives of the first-order sensitivities with respect to each variable $x_i$ are given by:
\begin{align}
    \bm{\underline{0}} = \bm{J} \left(\frac{\partial}{\partial x_i}\nabla_{\bm{x}} \bm{y} \right) + \left(\sum_{k} \left[\frac{\partial}{\partial y_k}\bm{J}\right] \frac{\partial y_k}{\partial x_i} \right) \nabla_{\bm{x}} \bm{y},
\end{align}
where $\bm{\underline{0}}$ is a zero matrix. Note that the second-order sensitivity of $y_k$ with respect to $x_i$ and $x_l$ is given by
\begin{align}
    \left[\frac{\partial}{\partial x_i}\nabla_{\bm{x}} \bm{y} \right]_{kl}= \frac{\partial^2 y_k}{\partial x_i \partial x_l}.
\end{align}

Thus, we can compute
\begin{align}
    \left[\frac{\partial}{\partial x_i}\nabla_{\bm{x}} \bm{y} \right]= - \bm{J}^{-1} \left(\sum_{k} \textcolor{black}{\left[\frac{\partial}{\partial y_k}\bm{J}\right]} \frac{\partial y_k}{\partial x_i} \right) \bm{J}^{-1}. \label{eq:sensitivity}
\end{align}

We next show how to compute the terms in~\eqref{eq:sensitivity}. Let diag$(\cdot)$ represent the diagonal matrix with the vector argument on the diagonal. We denote $\text{Re}(\cdot)$ and $\text{Im}(\cdot)$ as the real and imaginary part of the quantity $(\cdot)$, respectively. The asterisk $(*)$ represents the complex conjugate. We can now express the Jacobian matrix in terms of the derivative of complex power $S$ (i.e., $P + jQ$) with respect to the voltage magnitudes and angles:
\begin{equation}
\bm{J} = 
\begin{bmatrix}
    \text{Re}\left\{\nabla_{\bm{\theta}} \bm{S} \right\} & \text{Re}\left\{\nabla_{\bm{V}} \bm{S} \right\} \\
    \text{Im}\left\{\nabla_{\bm{\theta}} \bm{S} \right\} & \text{Im}\left\{\nabla_{\bm{V}} \bm{S} \right\}
    \label{eq:jocobian}
\end{bmatrix}.
\end{equation}

Let $\bm{Y}$ be the network admittance matrix. To compute the partial derivatives of the Jacobian matrix, i.e., $\frac{\partial}{\partial y_k}\bm{J}$ in~\eqref{eq:sensitivity}, we differente the entries of the Jacobian, as shown in~\cite{matpowerTechNote2}:
\begin{subequations}
\label{eq:jacobian}
\vspace{-0.25em}
\begin{align}
	\nabla_{\bm{\theta}} \bm{S} &= j \ \text{diag}\left(\bm{V}\right) \left(\text{diag}\left(\bm{I}\right) - \bm{Y} \text{diag}\left(\bm{V}\right)\right)^*, \label{eq:jacobian_t} \\
    \nabla_{\bm{V}} \bm{S}  &= \text{diag}\left({e}^{j\bm{\theta}}\right)\text{diag}\left(\bm{I}\right)^*  + \text{diag}(\bm{V})\left(\bm{Y}\text{diag}\left({e}^{j\bm{\theta}}\right)\right)^*, \label{eq:jacobian_v} 
\end{align}
\end{subequations}
where $\bm{I} = \bm{YV}$ is the bus current injection vector. We obtain:
\begin{subequations}
\label{eq:jacobian_derivative}
\vspace{-0.25em}
\begin{align}
	\frac{\partial}{\partial \theta_k}\left(\nabla_{\bm{\theta}} \bm{S}\right) &= \text{diag}\bm{\underline{V}_m} (-\text{diag}\left(\bm{I}\right)^* + \bm{Y} \text{diag}(\bm{V}))^* \notag \\
 & \ + \text{diag}(\bm{V}) (\text{diag}(\bm{Y} \bm{\underline{V}_k}) - \bm{Y} \text{diag}\bm{\underline{V}_k})^*, \label{eq:jacobian_t_tm} \\
    \frac{\partial}{\partial \theta_k}\left(\nabla_{\bm{V}} \bm{S}\right) &= j(-\text{diag}({e}^{j\bm{\theta}})(\text{diag}(\bm{Y} \bm{\underline{V}_k}))^* + \notag \\  & \quad \text{diag}(\underline{e}^{j\theta_k})\text{diag}(\bm{I})^*   
 - \text{diag}(\bm{V})(\bm{Y}\text{diag}(\underline{e}^{j\theta_k}))^* \notag \\  & \quad+ \text{diag}\bm{\underline{V}_k}(\bm{Y}\text{diag}({e}^{j\bm{\theta}}))^*), \label{eq:jacobian_v_tm} \\
  \frac{\partial}{\partial V_k}\left(\nabla_{\bm{\theta}} \bm{S}\right) &= j(\text{diag}(\underline{e}^{j\theta_k})(\text{diag}\left(\bm{I}\right) - \bm{Y} \text{diag}(\bm{V}))^*  \notag \\ 
 & \ + \text{diag}(\bm{V})(\text{diag}(\bm{Y} \underline{e}^{j\theta_k})^* - (\bm{Y}\text{diag}(\underline{e}^{j\theta_k}))^*), \label{eq:jacobian_t_vm} \\
 \frac{\partial}{\partial V_k}\left(\nabla_{\bm{V}} \bm{S}\right) &= \text{diag}\left({e}^{j\bm{\theta}}\right)\text{diag}\left(\bm{Y} \underline{e}^{j\theta_k}\right)^* \notag \\
 & \quad + \text{diag}\left(\underline{e}^{j\theta_k}\right)\left(\bm{Y}\text{diag}\left({e}^{j\bm{\theta}}\right)\right)^*, \label{eq:jacobian_v_vm} 
\end{align}
\end{subequations}
where $\bm{\underline{V}_k}$ is an all zero vector except for the $k^{\text{th}}$ entry which has the value $V_k$. Similarly, $\underline{e}^{j\theta_k}$ represents an all zero vector except for the $k^{\text{th}}$ entry which has the value $e^{j\theta_k}$. Let $\bm{\Gamma}_{(\cdot)}$ denote the derivative of the Jacobian matrix with respect to the subscript $(\cdot)$, i.e., $\frac{\partial}{\partial (\cdot)}\bm{J}$ in~\eqref{eq:sensitivity}. The derivative of the Jacobian matrix with respect to, e.g., the voltage angle $\theta_k$ is given by:
\begin{equation}
    \bm{\Gamma}_{\theta_k} = 
    \begin{bmatrix}
    \text{Re}\left\{\cfrac{\partial}{\partial \theta_k} \left(\nabla_{\bm{\theta}} \bm{S} \right)\right\} & \text{Re}\left\{\cfrac{\partial}{\partial \theta_k} \left(\nabla_{\bm{V}} \bm{S} \right)\right\} \\
    \text{Im}\left\{\cfrac{\partial}{\partial \theta_k} \left(\nabla_{\bm{\theta}} \bm{S} \right)\right\} & \text{Im}\left\{\cfrac{\partial}{\partial \theta_k} \left(\nabla_{\bm{V}} \bm{S} \right)\right\}
\end{bmatrix}.
\end{equation}
We can similarly compute the derivative of the Jacobian matrix with respect to a voltage magnitude.

Let $[\bm{J}^{-1}]_k$ represent the $k^{\text{th}}$ row of the inverse Jacobian matrix $\bm{J}^{-1}$. Then, the second-order sensitivity matrix $\bm{\Lambda}_{y_k}$ is:
\begin{equation}
\bm{\Lambda}_{y_k} = 
\begin{bmatrix}
    - [\bm{J}^{-1}]_k \left(\sum_{m} \bm{\Gamma}_{y_m} \frac{\partial y_m}{\partial x_1} \right) \bm{J}^{-1} \\
    - [\bm{J}^{-1}]_k \left(\sum_{m} \bm{\Gamma}_{y_m} \frac{\partial y_m}{\partial x_2} \right) \bm{J}^{-1} \\
    \vdots \\
    - [\bm{J}^{-1}]_k \left(\sum_{m} \bm{\Gamma}_{y_m} \frac{\partial y_m}{\partial x_{2N}} \right) \bm{J}^{-1} 
\end{bmatrix},
\end{equation}

\bibliographystyle{IEEEtran}
\bibliography{reference.bib}

\begin{thebibliography}{10}
\providecommand{\url}[1]{#1}
\csname url@samestyle\endcsname
\providecommand{\newblock}{\relax}
\providecommand{\bibinfo}[2]{#2}
\providecommand{\BIBentrySTDinterwordspacing}{\spaceskip=0pt\relax}
\providecommand{\BIBentryALTinterwordstretchfactor}{4}
\providecommand{\BIBentryALTinterwordspacing}{\spaceskip=\fontdimen2\font plus
\BIBentryALTinterwordstretchfactor\fontdimen3\font minus \fontdimen4\font\relax}
\providecommand{\BIBforeignlanguage}[2]{{%
\expandafter\ifx\csname l@#1\endcsname\relax
\typeout{** WARNING: IEEEtran.bst: No hyphenation pattern has been}%
\typeout{** loaded for the language `#1'. Using the pattern for}%
\typeout{** the default language instead.}%
\else
\language=\csname l@#1\endcsname
\fi
#2}}
\providecommand{\BIBdecl}{\relax}
\BIBdecl

\bibitem{stott2009}
B.~Stott, J.~Jardim, and O.~Alsa\c{c}, ``{DC} power flow revisited,'' \emph{IEEE Transactions on Power Systems}, vol.~24, no.~3, pp. 1290--1300, Aug. 2009.

\bibitem{baran1989}
M.~E. Baran and F.~F. Wu, ``Optimal capacitor placement on radial distribution systems,'' \emph{IEEE Transactions on Power Delivery}, vol.~4, no.~1, pp. 725--734, Jan. 1989.

\bibitem{taheri_molzahn-dcparam}
B.~Taheri and D.~K. Molzahn, ``Optimizing parameters of the {DC} power flow,'' \emph{Electric Power Systems Research}, vol. 235, no. 110719, October 2024, {\rm presented at the} \textit{23rd Power Systems Computation Conference (PSCC)}.

\bibitem{yang2018}
Z.~Yang, H.~Zhong, A.~Bose, T.~Zheng, Q.~Xia, and C.~Kang, ``A linearized {OPF} model with reactive power and voltage magnitude: A pathway to improve the {MW}-only {DC OPF},'' \emph{IEEE Transactions on Power Systems}, vol.~33, no.~2, pp. 1734--1745, Mar. 2018.

\bibitem{fan2021}
Z.~Fan, Z.~Yang, J.~Yu, K.~Xie, and G.~Yang, ``Minimize linearization error of power flow model based on optimal selection of variable space,'' \emph{IEEE Transactions on Power Systems}, vol.~36, no.~2, pp. 1130--1140, Mar. 2021.

\bibitem{BUASON2022}
P.~Buason, S.~Misra, and D.~K. Molzahn, ``A sample-based approach for computing conservative linear power flow approximations,'' \emph{Electric Power Systems Research}, vol. 212, p. 108579, 2022, {\rm presented at the} \emph{22nd Power Systems Computation Conference (PSCC)}.

\bibitem{fnt}
D.~K. Molzahn and I.~A. Hiskens, ``A survey of relaxations and approximations of the power flow equations,'' \emph{Foundations and Trends in Electric Energy Systems}, vol.~4, no. 1-2, pp. 1--221, Feb. 2019.

\bibitem{misra2018optimal}
S.~Misra, D.~K. Molzahn, and K.~Dvijotham, ``Optimal adaptive linearizations of the {AC} power flow equations,'' in \emph{20th Power Systems Computation Conference (PSCC)}, June 2018.

\bibitem{muhlpfordt2019optimal}
T.~M{\"u}hlpfordt, V.~Hagenmeyer, D.~K. Molzahn, and S.~Misra, ``Optimal adaptive power flow linearizations: Expected error minimization using polynomial chaos expansion,'' in \emph{IEEE Milan PowerTech}, 2019.

\bibitem{taheri_gupta_molzahn-opt_lindistflow}
B.~Taheri, R.~Gupta, and D.~K. Molzahn, ``Optimized {LinDistFlow} for high-fidelity power flow modeling of distribution networks,'' \emph{arXiv:2404.05125}, 2024.

\bibitem{buason2023datadriven}
P.~Buason, S.~Misra, S.~Talkington, and D.~K. Molzahn, ``A data-driven sensor placement approach for detecting voltage violations in distribution systems,'' \emph{Electric Power Systems Research}, vol. 232, no. 110387, 2024.

\bibitem{buason_misra_molzahn-cbla}
P.~Buason, S.~Misra, and D.~K. Molzahn, ``{Sample-Based Conservative Bias Linear Power Flow Approximations},'' in \emph{IEEE IAS Industrial and Commercial Power System Asia (IEEE I\&CPS Asia)}, July 2024.

\bibitem{CHEN2022108573}
J.~Chen and L.~A. Roald, ``A data-driven linearization approach to analyze the three-phase unbalance in active distribution systems,,'' \emph{Electric Power Systems Research}, vol. 211, p. 108573, 2022, {\rm presented at the} \emph{22nd Power Systems Computation Conference (PSCC)}.

\bibitem{liu2018}
Y.~Liu, N.~Zhang, Y.~Wang, J.~Yang, and C.~Kang, ``Data-driven power flow linearization: A regression approach,'' \emph{IEEE Transactions on Smart Grid}, vol.~10, no.~3, pp. 2569--2580, May 2019.

\bibitem{10202779}
M.~Jia and G.~Hug, ``Overview of data-driven power flow linearization,'' in \emph{IEEE Belgrade PowerTech}, 2023.

\bibitem{jia2023tutorial1}
M.~Jia, G.~Hug, N.~Zhang, Z.~Wang, and Y.~Wang, ``Tutorial on data-driven power flow linearization--{Part I}: Challenges and training algorithms,'' \emph{{\rm preprint available at \url{https://doi.org/10.3929/ethz-b-000606654}}}, 2023.

\bibitem{jia2023tutorial2}
------, ``Tutorial on data-driven power flow linearization--{Part II}: {S}upportive techniques and experiments,'' \emph{{\rm preprint available at \url{https://doi.org/10.3929/ethz-b-000606656}}}, 2023.

\bibitem{baker2021}
K.~Baker, ``Solutions of {DC OPF} are never {AC} feasible,'' in \emph{12th ACM International Conference on Future Energy Systems}, 2021, pp. 264--268.

\bibitem{dvijotham_molzahn-cdc2016}
K.~Dvijotham and D.~K. Molzahn, ``Error bounds on the {DC} power flow approximations: A convex relaxation approach,'' \emph{55th IEEE Conference on Decision and Control (CDC)}, December 2016.

\bibitem{castillo2016}
A.~Castillo, C.~Laird, C.~A. Silva-Monroy, J.-P. Watson, and R.~P. O’Neill, ``The unit commitment problem with {AC} optimal power flow constraints,'' \emph{IEEE Transactions on Power Systems}, vol.~31, no.~6, pp. 4853--4866, 2016.

\bibitem{bhusal2020}
N.~Bhusal, M.~Abdelmalak, M.~Kamruzzaman, and M.~Benidris, ``Power system resilience: Current practices, challenges, and future directions,'' \emph{IEEE Access}, vol.~8, pp. 18\,064--18\,086, 2020.

\bibitem{austgen2023comparisons}
B.~Austgen, E.~Kutanoglu, J.~J. Hasenbein, and S.~Santoso, ``Comparisons of two-stage models for flood mitigation of electrical substations,'' \emph{{\rm to appear in} INFORMS Journal on Computing}, 2024.

\bibitem{haag2024}
E.~Haag, N.~Rhodes, and L.~Roald, ``Long solution times or low solution quality: On trade-offs in choosing a power flow formulation for the optimal power shutoff problem,'' \emph{Electric Power Systems Research}, vol. 234, p. 110713, 2024, presented at \textit{23rd Power Systems Computation Conference (PSCC)}.

\bibitem{owen_aquino_talkington_molzahn-EVacuation_feeder}
A.~D. {Owen Aquino}, S.~Talkington, and D.~K. Molzahn, ``Managing vehicle charging during emergencies via conservative distribution system modeling,'' in \emph{Texas Power and Energy Conference (TPEC)}, Feb. 2024.

\bibitem{matpowerTechNote2}
R.~D. Zimmerman, ``{AC} power flows, generalized {OPF} costs and their derivatives using complex matrix notation,'' \emph{{\sc Matpower} Technical Note~2}, Feb. 2010.

\bibitem{lee2018}
J.-O. Lee, Y.-S. Kim, E.-S. Kim, and S.-I. Moon, ``Generation adjustment method based on bus-dependent participation factor,'' \emph{IEEE Transactions on Power Systems}, vol.~33, no.~2, pp. 1959--1969, 2018.

\bibitem{martin2016}
J.~A. Martin and I.~A. Hiskens, ``Generalized line loss relaxation in polar voltage coordinates,'' \emph{IEEE Transactions on Power Systems}, vol.~32, no.~3, pp. 1980--1189, 2017.

\bibitem{jabr2019}
R.~A. Jabr, ``High-order approximate power flow solutions and circular arithmetic applications,'' \emph{IEEE Transactions on Power Systems}, vol.~34, no.~6, pp. 5053--5062, 2019.

\bibitem{Cuyt1984}
A.~Cuyt, \emph{Abstract {Pad\'e} Approximants in Operator Theory}.\hskip 1em plus 0.5em minus 0.4em\relax Berlin, Heidelberg: Springer Berlin Heidelberg, 1984, pp. 1--58.

\bibitem{trias2018helm}
A.~Trias \emph{et~al.}, ``{HELM}: The holomorphic embedding load-flow method. {F}oundations and implementations,'' \emph{Foundations and Trends{\textregistered} in Electric Energy Systems}, vol.~3, no. 3-4, pp. 140--370, 2018.

\bibitem{BAKER196121}
G.~A. Baker and J.~Gammel, ``The {Padé} approximant,'' \emph{Journal of Mathematical Analysis and Applications}, vol.~2, no.~1, pp. 21--30, 1961.

\bibitem{cuyt1999}
A.~Cuyt, ``How well can the concept of {Pad{\'e}} approximant be generalized to the multivariate case?'' \emph{Journal of Computational and Applied Mathematics}, vol. 105, no. 1-2, pp. 25--50, 1999.

\bibitem{owen_aquino_roald_molzahn-ac_constraint_screening}
A.~D. {Owen Aquino}, L.~A. Roald, and D.~K. Molzahn, ``Identifying redundant constraints for {AC OPF}: The challenges of local solutions, relaxation tightness, and approximation inaccuracy,'' in \emph{53rd North American Power Symposium (NAPS)}, 2021.

\bibitem{gupta_buason_molzahn-fairness_pv_limit}
R.~Gupta, P.~Buason, and D.~K. Molzahn, ``Fairness-aware photovoltaic generation limits for voltage regulation in power distribution networks using conservative linear approximations,'' in \emph{8th Texas Power and Energy Conference (TPEC)}, February 2024.

\bibitem{zimmerman_matpower_2011}
R.~D. Zimmerman, C.~E. Murillo-Sánchez, and R.~J. Thomas, ``\mbox{MATPOWER}: Steady-state operations, planning, and analysis tools for power systems research and education,'' \emph{IEEE Transactions on Power Systems}, vol.~26, no.~1, pp. 12--19, Feb. 2011.

\bibitem{wogrin2020}
S.~Wogrin, S.~Pineda, and D.~A. Tejada-Arango, \emph{Applications of Bilevel Optimization in Energy and Electricity Markets}.\hskip 1em plus 0.5em minus 0.4em\relax Cham: Springer International Publishing, 2020, pp. 139--168.

\bibitem{Go2016}
R.~S. Go, F.~D. Munoz, and J.-P. Watson, ``Assessing the economic value of co-optimized grid-scale energy storage investments in supporting high renewable portfolio standards,'' \emph{Applied Energy}, vol. 183, no.~C, pp. 902--913, 2016.

\end{thebibliography}

\end{document}